\documentclass[11pt]{article}
\usepackage{psfrag,amsmath}
\usepackage{amssymb}
\usepackage{mathrsfs}
\usepackage{amsfonts}
\usepackage{amsmath}
\usepackage{mathrsfs,color}
\usepackage{epstopdf}
\usepackage{graphicx}
\usepackage{mathrsfs,amscd,amssymb,amsthm,amsmath,bm,url}

\makeatletter

\renewcommand{\@seccntformat}[1]{{\csname the#1\endcsname}{\normalsize .}\hspace{.5em}}
\makeatother

\def \[{\begin{equation}}
\def \]{\end{equation}}

\newtheorem{thm}{Theorem}[section]

\newtheorem{claim}{Claim}

\newtheorem{lem}[thm]{Lemma}
\newtheorem{cor}[thm]{Corollary}

\newtheorem{pb}{Problem}

\newenvironment{wst}
{\setlength{\leftmargini}{1.5\parindent}
 \begin{itemize}
 \setlength{\itemsep}{-1.1mm}}
{\end{itemize}}
\begin{document}

\setlength{\baselineskip}{15pt}

\begin{center}{\Large \bf Sharp bounds on the $A_{\alpha}$-index of graphs in terms of the independence number\footnote{
{\it Email addresses}: wtsun2018@sina.com (W. Sun),\ 289908600@qq.com (L. Yan),\ lscmath@ccnu.edu.cn (S. Li),\ xcli@uga.edu (X. Li).
}}
\vspace{4mm}

{\large Wanting Sun$^a$, \ \ Lixia Yan$^b$, \ \ Shuchao Li$^a$,\ \ Xuechao Li$^c$}\vspace{2mm}

$^a$ Faculty of Mathematics and Statistics,  Central China Normal
University, Wuhan 430079, P.R. China

$^b$ School of Information Eigeneering, Wuhan Business University, Wuhan 430056, P.R. China

$^c$  The University of Georgia, Athens, GA 30602, United States
\end{center}

\noindent {\bf Abstract}:\ Given a graph $G$, the adjacency matrix and degree diagonal matrix of $G$ are denoted by $A(G)$ and $D(G)$, respectively. In 2017,  Nikiforov \cite{0007} proposed the $A_{\alpha}$-matrix: $A_{\alpha}(G)=\alpha D(G)+(1-\alpha)A(G),$ where $\alpha\in [0, 1]$. The largest eigenvalue of this novel matrix is called the $A_\alpha$-index of $G$. In this paper, we characterize the graphs with  minimum $A_\alpha$-index among $n$-vertex graphs with independence number $i$ for $\alpha\in[0,1)$, where $i=1,\lfloor\frac{n}{2}\rfloor,\lceil\frac{n}{2}\rceil,{\lfloor\frac{n}{2}\rfloor+1},n-3,n-2,n-1,$ whereas for $i=2$ we consider the same problem for $\alpha\in [0,\frac{3}{4}{]}.$  Furthermore, we determine the unique graph (resp. tree) on $n$ vertices with given independence number having the maximum $A_\alpha$-index with $\alpha\in[0,1)$, whereas for the $n$-vertex bipartite graphs with given independence number, we characterize the unique graph having the maximum $A_\alpha$-index with $\alpha\in[\frac{1}{2},1).$

\vspace{2mm} \noindent{\it Keywords:}
Independence number; Spectral radius; $A_{\alpha}$-index; Bipartite graph
\vspace{2mm}


\section{\normalsize Introduction}\setcounter{equation}{0}
Throughout this paper, we assume $G$ is a simple (i.e., a finite, undirected, loopless and without multiple edges) connected graph, whose vertex set is $V_G=\{v_1,v_2,\ldots,v_n\}$ and edge set is $E_G$. The \textit{order} of $G$ is the number $n=|V_G|$ of its vertices and its \textit{size} is the number $|E_G|$ of its edges. The \textit{complement} of a graph $G$ is a graph $G^c$ with the same vertex set as $G$, in which any two distinct vertices are adjacent if and only if they are non-adjacent in $G$. Unless otherwise stated, we follow the traditional notation and terminology (see \cite{0009}).

The {\textit{degree} $d_G(u)$} of a vertex $u$ (in a graph $G$) is the number of edges incident with it. We say that two vertices $i$ and $j$ are \textit{adjacent} (or \textit{neighbours}) if they are joined by an edge. The \textit{set of neighbours} of a vertex $u$ is denoted by $N_G(u)$. Then $\Delta(G):=\max_{v\in V_G}d_G(v)$ is the maximum degree of $G$.  An edge is a \textit{pendant edge} if one of its ends has degree $1.$ A vertex of degree at least $3$ is a \textit{major vertex}. A path, a cycle, a complete bipartite graph and a complete graph on $n$ vertices are denoted by $P_n$, $C_n$, $K_{a,n-a}$ and $K_n$, respectively.

An \textit{internal path} of $G$ is a walk $v_0v_1\ldots v_s\, (s \geqslant 1)$ such that the vertices $v_0, v_1, \ldots, v_s$ are distinct, $d_G(v_0)\geqslant 3, d_G(v_s)\geqslant 3$ and $d_G(v_i)=2$, whenever $1\leqslant i\leqslant {s-1}$. For two disjoint graphs $G$ and $H$, the \textit{join} $G\vee H$ is the graph obtained by joining every vertex of $G$ with every vertex of
$H$.

An \textit{independent set} in a graph is a set of pairwise non-adjacent vertices, and the number of vertices
in the largest independent set is called the \textit{independence number}, denoted by $i(G)$. Let $\mathcal{G}_{n,i}$ be the set of connected graphs with fixed order $n$ and independence number $i$.  A \textit{matching} in a graph $G$ is a set of edges without common vertices. The matching number $\mu(G)$ is the maximum size of a matching in~$G$. 

Given a graph $G$ of order $n$, its \textit{adjacency matrix} $A(G)$ is an $n\times n$ $0$-$1$ matrix whose $(i,j)$-entry is $1$ if and only if $v_i$ is adjacent to $v_j$ in $G$. We say that the largest eigenvalue of $A(G)$ is the \textit{spectral radius} of $G.$ Let $D(G)={\rm diag}(d_1,\ldots, d_n)$ be the diagonal matrix of vertex degrees in a graph $G$. Then the matrix $Q(G)=D(G)+A(G)$ is called the \textit{signless Laplacian matrix} (see \cite{RD}). 
Quite recently, Nikiforov \cite{0007} proposed the \textit{$A_{\alpha}$-matrix} as a convex combination of $D(G)$ and $A(G)$, that is,
$$
  A_{\alpha}(G)=\alpha D(G)+(1-\alpha)A(G),\ \ \ 0\leqslant \alpha\leqslant 1.
$$
Notice that $A_{\alpha}(G)$ is real and symmetric. Hence its real eigenvalues can be given in non-increasing order as $\lambda_1(G)\geqslant  \cdots \geqslant \lambda_n(G)$. 
For short, the $A_{\alpha}$-spectral radius (i.e., the largest eigenvalue of $A_{\alpha}(G)$), denoted by $\lambda_\alpha(G)=\lambda_1(G),$ is called the \textit{$A_{\alpha}$-index} of $G.$ Notice that
$$
  A(G)=A_0(G),\ \ \ Q(G)=2A_{\frac{1}{2}}(G) \ \ \ \text{and}\ \ \ D(G)=A_1(G).
$$

Brualdi and Solheid \cite{0002} posed the following general problem, which became one of the classic problems in the theory of graph spectra:
\begin{pb}\label{pb1}
  Given a set of graphs $\mathcal{G}$, establish an upper or a lower bound for the spectral radius of graphs in $\mathcal{G},$ and characterize the graphs in which the maximum or the minimum spectral radius is achieved.
\end{pb}
For the maximization part of Problem \ref{pb1}, one may be referred to representative achievements \cite{002,009,008,013,005} and the references cited in.
For the minimization part of Problem \ref{pb1}, one may consult \cite{0010,0013,Xu}. In addition, similar problems for the signless Laplacian spectral radius also have been studied extensively. For more advances, one may be referred to \cite{0006,001,006,007,0004,016} and references therein.

According to Problem \ref{pb1}, it is natural to study the problem: Given a class of graphs $\mathcal{G}$ and $\alpha\in[0,1)$, find $\min\{\lambda_\alpha(G): G\in \mathcal{G}\}$ and $\max\{\lambda_\alpha(G): G\in \mathcal{G}\}$, and characterize the graphs achieving the minimum or maximum values.

Huang, Lin and Xue \cite{Huang} studied the Nordhaus-Gaddum type bounds for the $A_{\alpha}$-index of graphs, and they determined the graphs with maximum $A_{\alpha}$-index among all graphs with given order and size. Li, Chen and Meng \cite{Li} determined the unique tree (resp. the unicyclic graphs) with given degree sequence having the maximum $A_{\alpha}$-index. Lin, Huang and Xue \cite{0001} characterized the extremal graph with maximum $A_{\alpha}$-index with fixed order and number of cut vertices, and the extremal tree which attains the maximum $A_{\alpha}$-index with fixed order and matching number.  For more advances on the $A_{\alpha}$-index and some attractive feature of $A_{\alpha}$-matrices, one may be referred to \cite{Chen,LS3,LS2,LS1,LI2019,LI2020,LY,Ni,0005,RJ,RJ1,Wang,Xu2020} and the references cited in.

Very recently, Nikiforov et al. \cite{Nik} given several bounds on the $A_{\alpha}$-index of graphs, and they determined the unique tree with maximum (resp. minimum) $A_{\alpha}$-index among trees. Xue et al. \cite{Xue} determined the graphs with maximum (resp. minimum) $A_{\alpha}$-index among all connected graphs with given diameter (resp. clique number). Along this line and motivated by \cite{0004,Xu}, it is natural and interesting to study  graphs with given independence number having the maximum or minimum $A_{\alpha}$-index.

Let $T_{a,b,c}$ denote a $T$-shaped tree defined as a tree with a single vertex $u$ of degree $3$ such that $T_{a,b,c}-u=P_a\cup P_b\cup P_c\, (1\leqslant a\leqslant b\leqslant c).$ In addition, if $a=b=1,$ the tree $T_{1,1,n-3}$ is also called a \textit{snake}.
Our first main result characterizes the $n$-vertex graphs with independence number $i$ having the minimum $A_{\alpha}$-index, where $i=1,\left\lfloor\frac{n}{2}\right\rfloor,\lceil\frac{n}{2}\rceil,{\lfloor\frac{n}{2}\rfloor+1},n-1.$
\begin{thm}\label{thm1.1}
Let $G$ be in $\mathcal{G}_{n,i}$ with minimum $A_{\alpha}$-index and $\alpha\in[0,1).$
\begin{wst}
\item[{\rm (i)}] If $i=1,$ then $G\cong K_n;$
\item[{\rm (ii)}] If $i=\left\lfloor\frac{n}{2}\right\rfloor,$ then $G\cong P_n$ for even $n$ and $G\cong C_n$ for odd $n;$
\item[{\rm (iii)}] If $i=\left\lceil\frac{n}{2}\right\rceil,$ then $G\cong P_n;$
\item[{\rm (iv)}] If $i={\lfloor\frac{n}{2}\rfloor+1}$ with $n\geqslant 9,$ then {$G\cong P_n$ for odd $n$ and $G\cong T_{1,1,n-3}$ for even $n;$}
\item[{\rm (v)}]If $i=n-1,$ then $G\cong K_{1,n-1}.$
\end{wst}
\end{thm}

\begin{figure}[!ht]
\begin{center}
  \psfrag{a}{$s_1$}\psfrag{b}{$t_1$}\psfrag{c}{$s_2$}\psfrag{d}{$t_2$}
  \psfrag{1}{$v_1$}\psfrag{2}{$v_2$}\psfrag{3}{$v_3$}
   \psfrag{s}{$G_1(s_1,t_1)$}\psfrag{t}{$G_2(s_2,t_2)$}
  \includegraphics[width=65mm]{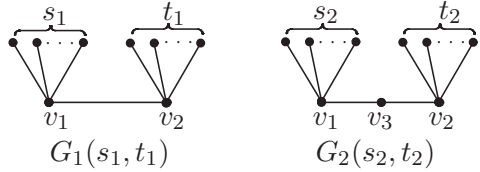}\\
  \caption{Graphs $G_1(s_1,t_1)$ and $G_2(s_2,t_2)$.}\label{f2}
  \end{center}
\end{figure}

Let $G_1(s_1,t_1)$ be a tree of order $n$ obtained from an edge $v_1v_2$ by attaching $s_1$ (resp. $t_1$) pendant edges to $v_1$ (resp.  $v_2$), where {$\min\{s_1,t_1\}\geqslant 1$ and} $s_1+t_1=n-2.$ Let $G_2(s_2,t_2)$ be a graph of order $n$ obtained from $G_1(s_2,t_2)$ by subdividing the edge $v_1v_2,$ where {$\max\{s_2,t_2\}\geqslant 1$ and} $s_2+t_2=n-3$ (see Fig. \ref{f2}).
\begin{figure}[!ht]
\begin{center}
  \psfrag{a}{$s_1$}\psfrag{b}{$t_1$}\psfrag{e}{$k_1$}\psfrag{c}{$s_2$}\psfrag{d}{$t_2$}\psfrag{o}{$k_2$}\psfrag{h}{$s_3$}\psfrag{i}{$t_3$}\psfrag{l}{$k_3$}\psfrag{x}{$s_4$}\psfrag{y}{$t_4$}\psfrag{z}{$k_4$}
  \psfrag{1}{$v_1$}\psfrag{2}{$v_2$}\psfrag{3}{$v_3$}\psfrag{4}{$v_4$}\psfrag{5}{$v_5$}
   \psfrag{s}{$H_1(s_1, t_1, k_1)$}\psfrag{t}{$H_2(s_2, t_2, k_2)$}\psfrag{r}{$H_3(s_3, t_3, k_3)$}\psfrag{k}{$H_4(s_4, t_4, k_4)$}
  \includegraphics[width=150mm]{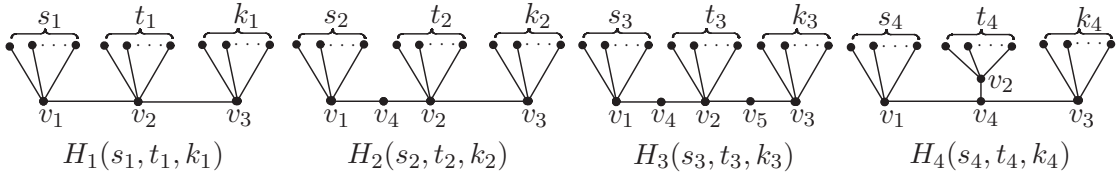}\\
  \caption{Graphs $H_1(s_1, t_1, k_1), H_2(s_2, t_2, k_2), H_3(s_3, t_3, k_3),$ and $H_4(s_4, t_4, k_4)$.}\label{f03}
  \end{center}
\end{figure}

\begin{figure}[!ht]
\begin{center}
  \psfrag{a}{\small $K_s$}\psfrag{b}{\small $K_t$}\psfrag{c}{\small $u$} \psfrag{d}{\small $v$}
  \includegraphics[width=50mm]{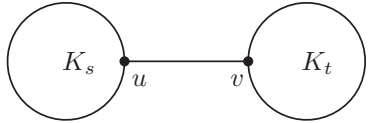}\\
  \caption{The graph $F(s,t)$.}\label{f3}
  \end{center}
\end{figure}

Let $H_1(s_1, t_1, k_1)$ be a tree obtained from a path $v_1v_2v_3$ by attaching $s_1$ (resp. $t_1$ and $k_1$) pendant edges to $v_1$ (resp. $v_2$ and $v_3$), where {$\min\{s_1,t_1,k_1\}\geqslant 1$ and} $s_1+t_1+k_1=n-3.$ Let $H_2(s_2, t_2, k_2)$ be a tree obtained from $H_1(s_2, t_2, k_2)$ by subdividing the edge $v_1v_2,$ where {$\max\{s_2,t_2\}\geqslant 1,$ $k_2\geqslant 1$ and} $s_2+t_2+k_2=n-4.$ Let $H_3(s_3, t_3, k_3)$ denote a tree obtained from $H_2(s_3, t_3, k_3)$ by subdividing the edge $v_2v_3,$ where {$\max\{s_3,t_3,k_3\}\geqslant 1$ and} $s_3+t_3+k_3=n-5.$ Let $H_4(s_4, t_4, k_4)=H_2(s_4, t_4, k_4)-v_2v_3+v_3v_4,$ where {at least two of $s_4,t_4,k_4$ are positive and} $s_4+t_4+k_4= n-4$ (see Fig. \ref{f03}). 
The next two main results characterize the $n$-vertex graphs with independence number $n-2$ or $n-3$ having the minimum $A_\alpha$-index, respectively.
\begin{thm}\label{thm3.1}
Let $G$ be in $\mathcal{G}_{n,n-2}$ with $n\geqslant 5$. If $\alpha\in [0,1),$ then
$$
  \lambda_\alpha(G)\geqslant \lambda_\alpha\left(G_2\left(\left\lceil\frac{n-3}{2}\right\rceil,\left\lfloor\frac{n-3}{2}\right\rfloor\right)\right)
$$
with equality if and only if $G\cong G_2\left(\left\lceil\frac{n-3}{2}\right\rceil,\left\lfloor\frac{n-3}{2}\right\rfloor\right).$
\end{thm}
\begin{thm}\label{thm3.2}
Let $G$ be in $\mathcal{G}_{n,n-3}$ with $n\geqslant 9$ having the minimum $A_\alpha$-index and let $\alpha\in [0,1).$ Then $G\cong H_3(s,t,k)$ for some integers $s,t,k$ {with $|s-k|\leqslant 1$ and $s+t+k=n-5.$} Furthermore,
\begin{wst}
\item[{\rm (i)}] if $\alpha\in [\frac{1}{2},1)$ and $n=3m,$ then  $G\cong H_3(m-1,m-3,m-1);$
\item[{\rm (ii)}] if $\alpha\in [\frac{1}{2},1)$ and $n=3m+1,$ then  $G\cong H_3(m-1,m-2,m-1);$
\item[{\rm (iii)}] if $\alpha\in [\frac{11}{20},1)$ and $n=3m+2,$ then $G\cong H_3(m,m-3,m).$
\end{wst}
\end{thm}
Denote $F(s,t):=(K_{s,t}-e)^c$ (see Fig. \ref{f3}). Our fourth main result identifies the unique graph among $\mathcal{G}_{n,2}$ having the minimum $A_\alpha$-index.
\begin{thm}\label{thm5.1}
Let $G$ be in $\mathcal{G}_{n,2}$ with $n\geqslant 11$ having the minimum $A_\alpha$-index. If $\alpha\in [0,\frac{3}{4}{]},$ then $G\cong F(\left\lceil\frac{n}{2}\right\rceil,\left\lfloor\frac{n}{2}\right\rfloor).$
\end{thm}
For $n\geqslant 2k+1,$ let $S_{n,n-k}^*$ be the graph obtained from the star $K_{1,n-k}$ by attaching a pendant edge to each of certain $k-1$ non-central vertices of $K_{1,n-k}.$ The last two main results characterize the unique $n$-vertex tree (resp. bipartite graph, simple graph) with independence number $i$ which attains the maximum $A_\alpha$-index among trees (resp. bipartite graphs, simple graphs).
\begin{thm}\label{thm1.5}
Let $G$ be an $n$-vertex bipartite graph with independence number $i$, where $\left\lceil\frac{n}{2}\right\rceil\leqslant i\leqslant n-1.$ If $\alpha\in [\frac{1}{2},1),$ then
$
  \lambda_\alpha(G)\leqslant \lambda_\alpha(K_{i,n-i})
$
with equality if and only if $G\cong K_{i,n-i}.$ In particular, assume that $G$ is an $n$-vertex tree with $n\geqslant 4.$
If $\alpha\in[0,1),$ then
$
  \lambda_\alpha(G)\leqslant \lambda_\alpha(S_{n,i}^*)
$
with equality if and only if $G\cong S_{n,i}^*.$
\end{thm}
\begin{thm}\label{thm1.6}
Let $G$ be an $n$-vertex graph with independence number $i$. If $\alpha\in[0,1),$ then
$
  \lambda_\alpha(G)\leqslant \lambda_\alpha(K_i^c\vee K_{n-i})
$
with equality if and only if $G\cong K^c_i\vee K_{n-i}.$
\end{thm}
The remainder of our paper is organized as follows: In Section 2, we review some preliminaries that will be needed in the sequel. In Section 3, we give the proof for Theorem \ref{thm1.1}. In Section 4, we give the proofs of Theorems \ref{thm3.1} and \ref{thm3.2}. In Section 5, we give the proof of Theorem \ref{thm5.1}. In Section 6, we present the proofs of Theorems \ref{thm1.5} and \ref{thm1.6}. In the last section, we give some conclusion remarks and some further research problems.
\section{\normalsize Preliminaries}
\begin{figure}[!ht]
\begin{center}
  \includegraphics[width=100mm]{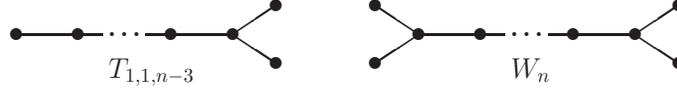}\\
  \caption{Graphs $T_{1,1,n-3}$ and $W_n$.}\label{f4}
  \end{center}
\end{figure}

\setcounter{equation}{0}
In this section, we present some preliminary results, which will be used in the subsequent sections. 
\begin{lem}[\cite{0001}]\label{lem1.1}
Let $T$ be a tree with order $n\geqslant 4$ and matching number $\mu\ \left(1\leqslant \mu\leqslant \left\lfloor\frac{n}{2}\right\rfloor\right)$. If $\alpha\in[0,1],$ then
$$
  \lambda_\alpha(T)\leqslant \lambda_\alpha(S_{n,n-\mu}^*)
$$
with equality if and only if $T\cong S_{n,n-\mu}^*.$
\end{lem}
\begin{lem}[\cite{0007}]\label{lem3.2}
Let $G$ be a connected graph and $H$ be a proper subgraph of $G.$ Then $\lambda_\alpha(H)<\lambda_\alpha(G)$ for $\alpha\in[0,1).$
\end{lem}
Note that $A_\alpha(G)$ is a non-negative and irreducible matrix if $G$ is connected. By the famous Perron-Frobenius Theorem, there exists a positive eigenvector, say ${\bf x}=(x_1,x_2,\ldots,x_n)^T,$ of $A_\alpha(G)$ corresponding to $\lambda_\alpha(G)$. For a vector ${\bf x}$ on vertices of $G,$ we denote by $x_{v_i}$ or $x_i$ the entry of ${\bf x}$ at $v_i\in V_G.$ 
\begin{lem}[\cite{Xue}]\label{lem1.2}
Let $G$ be a graph on $n$ vertices with $\alpha\in[0,1).$ For $v,u\in V_G,$ assume $S\subset N_G(v)\setminus (N_G(u)\cup \{u\}).$ Let $G'=G-\{vw:w\in S\}+\{uw:w\in S\}.$ Let ${\bf x}:=(x_1,x_2,\ldots,x_n)^T$ be a positive eigenvector of $A_\alpha(G)$ corresponding to $\lambda_\alpha(G).$ If $S\neq \emptyset$ and $x_u\geqslant x_v$, then $\lambda_\alpha(G')>\lambda_\alpha(G).$
\end{lem}
\begin{lem}[\cite{Nik}]\label{lem2.1}
Let $G$ be a connected graph on $n$ vertices and $\alpha\in[0,1)$. Then
$$
  \lambda_\alpha(G)\geqslant \lambda_\alpha(P_n)
$$
with equality if and only if $G\cong P_n.$
\end{lem}
\begin{lem}[\cite{0007}]\label{lem3.7}
Let $G$ be a graph with maximum degree $\Delta$ and let $\alpha\in[0,1)$. Then
\[\label{eq1.1}
  \lambda_\alpha(G)\geqslant \frac{1}{2}\left(\alpha(\Delta+1)+\sqrt{\alpha^2(\Delta+1)^2+4\Delta(1-2\alpha)}\right)
\]
with equality if and only if $G\cong K_{1,\Delta}.$ In particular,
\begin{equation}\label{eq:1.1}
    \lambda_\alpha(G)\geqslant
    \left\{
    \begin{aligned}
        &\alpha(\Delta+1),\ \ \ \ \textrm{if $\alpha\in [0,\frac{1}{2}];$}\\
        &\alpha\Delta+\frac{(1-\alpha)^2}{\alpha},\ \ \textrm{if $\alpha\in [\frac{1}{2},1).$}\\
    \end{aligned}
    \right.
\end{equation}
\end{lem}
{It is straightforward to check that the bound in \eqref{eq1.1} (see, \cite[Proposition 12]{0007}) is better than that of \eqref{eq:1.1} (see, \cite[Corollary 13]{0007}) if $\alpha\neq \frac{1}{2},$ and they are equal if $\alpha= \frac{1}{2}.$ Hence the equality in \eqref{eq:1.1} holds if and only if $\alpha= \frac{1}{2}$ and $G\cong K_{1,\Delta}.$}
\begin{lem}[\cite{0007}]\label{lem3.8}
Let $G$ be a graph.  If $\alpha\in [0,1),$ then
$$
\frac{2|E_G|}{|V_G|}\leqslant \lambda_\alpha(G)\leqslant \max_{uv\in E_G}\{\alpha d_G(u)+(1-\alpha)d_G(v)\},
$$
{the first equality holds if and only if $G$ is regular.}
\end{lem}
\begin{lem}[\cite{Li}]\label{lem3.1}
Let $G$ be a connected graph and $uv$ be some edge on an internal path of $G$. Let $G_{uv}$ denote the graph obtained from $G$ by subdividing the edge $uv$ into edges $uw$ and $wv$. Then $\lambda_\alpha(G_{uv})\leqslant \lambda_\alpha(G)$ {for $\alpha\in[0,1).$}
\end{lem}
Let $H$ be a connected graph with $v_0\in V_H$, and let $P_k=v_1 v_2 \ldots v_k.$ Denote by $G_H(k)$ the graph obtained from $H$ and $P_k$ by joining $v_0$ and $v_1$ with an edge. Then let $G_H(k,s)$ be the graph obtained from $G_H(k)$ and a path $P_s$ by connecting $v_0$ and an end-vertex of $P_s$ with an edge.
\begin{lem}[\cite{Xue}]\label{lem3.3}
Let $G_H(k)$ be the graph defined above and ${\bf x}$ be a positive eigenvector of $A_\alpha(G_H(k))$ corresponding to $\lambda_\alpha(G_H(k))$. If $\alpha\in [0,1)$ and $\lambda_\alpha(G_H(k))>2$, then $x_{v_i}>x_{v_{i+1}}$ for each $i\in \{0,1,\ldots,\linebreak k-1\}$.
\end{lem}
\begin{lem}[\cite{Xue}]\label{lem3.10}
Let $G_H(k, s)$ be the graph defined above with $k\geqslant s+2.$ If $\alpha\in [0,1)$ and $\lambda_\alpha(G_H(k, s))\geqslant 2,$ then $\lambda_\alpha(G_H(k,s))<\lambda_\alpha(G_H(k-1,s+1)).$
\end{lem}
An $n$-vertex tree $W_n$ containing just two vertices, say $u$ and $v$, of degree $3$ is called a \textit{double snake} if deleting $u,v$ from $W_n$ yields 4 isolated vertices and a path $P_{n-6}$. The snake $T_{1,1,n-3}$ and the double snake $W_n$ are depicted in Fig.~\ref{f4}.

\begin{lem}[\cite{Wang1}]\label{lem1.3}
Let $G$ be a connected graph of order $n$. Then
\begin{wst}
\item[{\rm (i)}] $\lambda_\alpha(G)<2$ if and only if $G$ is one of the following graphs:
\begin{wst}
  \item[{\rm (a)}] $P_n\,(n\geqslant 1)$ for $\alpha\in [0,1);$
  \item[{\rm (b)}] $T_{1,1,n-3}\,(n\geqslant 4)$ for $\alpha\in [0,s_1);$
  \item[{\rm (c)}] $T_{1,2,2}$ for $\alpha\in [0,s_2),$ $T_{1,2,3}$ for $\alpha\in [0, s_3)$ and $T_{1,2,4}$ for $\alpha\in [0, s_4)$.
\end{wst}
\item[{\rm (ii)}] $\lambda_\alpha(G)=2$ if and only if $G$ is one of the following graphs:
\begin{wst}
 \item[{\rm (a)}] $C_n\,(n\geqslant 3)$ for $\alpha\in [0,1];$
 \item[{\rm (b)}] $P_n\,(n\geqslant 3)$ for $\alpha=1;$
 \item[{\rm (c)}] $W_n\,(n\geqslant 6)$ for $\alpha=0;$
 \item[{\rm (d)}] $T_{1,1,n-3}$ for $\alpha=s_1,$
 \item[{\rm (e)}] $T_{1,2,2}$ for $\alpha=s_2$, $T_{1,2,3}$ for $\alpha=s_3,$ $T_{1,2,4}$ for $\alpha=s_4;$
 \item[{\rm (f)}] $T_{1,3,3}$ for $\alpha= 0,$ $T_{1,2,5}$ for $\alpha=0,$ $K_{1,4}$ for $\alpha=0,$ and $T_{2,2,2}$ for $\alpha=0,$
\end{wst}
\end{wst}
where $s_1=\frac{4}{n+1+\sqrt{(n+1)^2-16}}$, $s_2=0.2192+$ is the solution of $2\alpha^3-11\alpha^2+16\alpha-3=0,$ $s_3=0.1206+$ is the solution of $\alpha^3-6\alpha^2+9\alpha-1=0$ and $s_4=0.0517+$ is the solution of $2\alpha^3-13\alpha^2+20\alpha-1=0.$
\end{lem}
Let $G$ be a graph with $u,v\in V_G.$ We say $u$ and $v$ are \textit{equivalent} in $G$ if there exists an automorphism $\varphi:G\rightarrow G$ such that $\varphi(u)=v.$
\begin{lem}[\cite{0007}]\label{lem3.09}
Let $G$ be a connected graph on $n$ vertices, and let $u$ and $v$ be equivalent vertices in $G.$ If ${\bf x}$ is a positive eigenvector of $A_\alpha(G)$ corresponding to $\lambda_\alpha(G)$, then $x_u=x_v.$
\end{lem}
Let $h(x)$ be a real polynomial in the variable $x$, and let $\lambda(h)$ be the largest real root of the equation $h(x)=0.$ Then the next lemma follows immediately.
\begin{lem}\label{lem3.9}
Let $g(x)$ and $h(x)$ be two monic polynomials with real roots. If $h(x)<g(x)$ for $x= \lambda(g),$ then $\lambda(h)>\lambda(g).$
\end{lem}
Let $M$ be a real symmetric matrix, whose columns and rows are indexed by $U=\{1,2,\ldots,n\}.$ Assume that $\pi:=U_1\cup U_2\cup \cdots \cup U_t$ is a partition of $U$. Then $M$ can be partitioned based on $\pi$ as
\begin{equation*}
    M=\left(
        \begin{array}{ccc}
          M_{11} & \cdots & M_{1t} \\
          \vdots & \ddots & \vdots \\
          M_{t1} & \cdots & M_{tt} \\
        \end{array}
      \right),
\end{equation*}
where $M_{ij}$ denotes the submatrix of $M$, indexed by the rows and columns of $U_i$ and $U_j$, respectively. Let $\pi_{ij}$ be the average row sum of $M_{ij}$ for $1\leqslant i,j \leqslant t$. The matrix $M_{\pi}=(\pi_{ij})$ is called the \textit{quotient matrix} of $M$. Moreover, we call $\pi$ an \textit{equitable partition} if the row sum of $M_{ij}$ is constant for $1\leqslant i,j \leqslant t$.
\begin{lem}[\cite{YLH}]\label{lem4.01}
Let $M$ be a real symmetric nonnegative matrix with an equitable partition $\pi$ and let $M_{\pi}$ be the corresponding quotient matrix. Then the largest eigenvalue of $M$ is equal to that of $M_{\pi}$.
\end{lem}
\begin{lem}[\cite{Xu}]\label{lem4.1}
Let $G$ be a connected graph on $n$ vertices. If the independence number of $G$ is $2$, then $|E_{G^c}|\leqslant \lfloor\frac{n^2}{4}\rfloor-1.$
\end{lem}

\section{\normalsize Proof of Theorem \ref{thm1.1}}
In this section, we give the proof of Theorem \ref{thm1.1}, which identifies the graphs with  minimum $A_{\alpha}$-index in $\mathcal{G}_{n,i}$ for  $i=1,\left\lfloor\frac{n}{2}\right\rfloor,\lceil\frac{n}{2}\rceil,{\lfloor\frac{n}{2}\rfloor+1},n-1.$ Further on the following lemma is needed.

\begin{lem}\label{lem1.4}
Let $G$ be a connected graph of order $n\,(n\geqslant 9).$ If $\alpha\in[0,1)$ and $G\not\in\{P_n,C_n\},$ then
$$
\lambda_\alpha(G)\geqslant \lambda_\alpha(T_{1,1,n-3}).
$$
The equality holds if and only if $G\cong T_{1,1,n-3}.$
\end{lem}
\begin{proof}
Let $G\not\in\{P_n,C_n\}$ be a connected graph with  minimum $A_\alpha$-index. By Lemma \ref{lem3.2}, we know that $G$ is a tree. If $\alpha\in[0,s_1],$ then by Lemma \ref{lem1.3} one has $G\cong T_{1,1,n-3},$ as desired. In what follows, we consider the case $\alpha\in(s_1,1).$

Since $G\not\cong P_n,$ $G$ contains at least one major vertex.  Choose a major vertex of $G$, say $u$, such that all but at most one graphs in $G[V_{G_j}\cup \{u\}]$ are paths, where $G_j$ denotes the connected component of $G-u$ for $1\leqslant j\leqslant d_G(u).$  Then assume, without loss of generality, that $G-u\cong G_{1}\cup P_{a_1}\cup \cdots\cup P_{a_s}$, here $s=d_G(u)-1\geqslant 2$ and $a_1\geqslant \cdots\geqslant a_s\geqslant 1.$ {Denote by $w$ and $v$ the pendant vertex of $G$ in $V_{P_{a_1}}$ and the vertex in $N_G(u)\cap V_{P_{a_s}},$ respectively. Let $H=G-uv+wv.$} 
Assume that ${\bf x}$ is a positive  eigenvector of $A_\alpha(H)$ corresponding to $\lambda_\alpha(H).$

Now, we show that there exists exactly one major vertex in $G.$ Otherwise, $G$ contains at least two major vertices. Then $H$ contains at least one major vertex. It follows that $H\not\in \{P_n,C_n\}.$ Notice that $n\geqslant 9$ and $\alpha\in (s_1,1).$ In view of Lemma \ref{lem1.3}, one has $\lambda_\alpha(H)>2.$ Together with Lemma \ref{lem3.3}, we have $x_{u}\geqslant x_{w}.$ By Lemma \ref{lem1.2}, we obtain $\lambda_\alpha(H)<\lambda_\alpha(G),$ which contradicts the minimality of $G$. Hence $G$ contains exactly one major vertex $u.$

If $d_G(u)\geqslant 4,$ then $d_H(u)\geqslant 3.$ It follows from Lemma~\ref{lem1.3} that $\lambda_\alpha(H)>2.$ Together with Lemmas~\ref{lem1.2} and \ref{lem3.3}, we obtain $\lambda_\alpha(H)<\lambda_\alpha(G),$ a contradiction. Hence $d_G(u)\leqslant 3,$ that is, $G\cong T_{a,b,c}$ with $a\leqslant b \leqslant c$ and $a+b+c=n-1.$ By Lemma~\ref{lem1.3}, one has $\lambda_\alpha(T_{a,b,c})>2$ {if $\alpha\in(s_1,1)$ and $n\geqslant 9.$} In view of Lemma \ref{lem3.10}, we have $G\cong T_{1,1,n-3}$. This completes the proof.
\end{proof}
Now, we are ready to prove Theorem \ref{thm1.1}.
\begin{proof}[\bf Proof of Theorem \ref{thm1.1}]\
It is routine to check that $\mathcal{G}_{n,1}=\{K_n\}$ and $\mathcal{G}_{n,n-1}=\{K_{1,n-1}\}.$ Hence (i) and (v) hold. Note that $P_n\in \mathcal{G}_{n,\lceil\frac{n}{2}\rceil}$ and $C_n\in \mathcal{G}_{n,\lfloor\frac{n}{2}\rfloor}.$ Together with Lemmas~\ref{lem2.1} and \ref{lem1.3}, we obtain (ii)-(iii) immediately. Furthermore, it is easy to see that $T_{1,1,n-3}\in \mathcal{G}_{n,{\lfloor\frac{n}{2}\rfloor+1}}.$ Together with Lemmas~\ref{lem2.1} and \ref{lem1.4}, then (iv) follows immediately.
\end{proof}

\section{\normalsize Proofs of Theorems \ref{thm3.1} and \ref{thm3.2}}
In this section, we give the proofs of Theorems \ref{thm3.1} and \ref{thm3.2}, which characterize graphs with  minimum $A_\alpha$-index among $\mathcal{G}_{n,n-2}$ and $\mathcal{G}_{n,n-3},$ respectively.

In order to prove Theorem \ref{thm3.1}, we need the following two lemmas.
\begin{lem}\label{lem3.4}
Let $G$ be a connected graph with $u,v\in V_G$ and $G-u\cong G-v.$ Assume $G_{s,t}\,(s\geqslant t\geqslant 1)$ is a graph obtained from $G$ by attaching $s$ new pendant edges $\{uw_1, uw_2,\ldots, uw_s\}$ at $u$ and $t$ new pendant edges $\{v w_{s+1}, v w_{s+2},\ldots,v w_{s+t}\}$ at $v.$ If $\alpha\in[0,1),$ then $\lambda_\alpha(G_{s,t})<\lambda_\alpha(G_{s+1,t-1}).$
\end{lem}
\begin{proof}
Let ${\bf x}$ be a positive eigenvector of $A_\alpha(G_{s,t})$ corresponding to $\lambda_\alpha(G_{s,t}).$ Firstly, we prove that $x_u\geqslant x_v.$ Suppose to the contrary that $x_u<x_v.$ Then $s\geqslant t+1.$ Otherwise, $s=t$ and so $u,\,v$ are equivalent vertices in $G_{s,t}.$ By Lemma \ref{lem3.09}, one has $x_u=x_v,$ a contradiction. Let $G_{t,s}=G_{s,t}-\{uw_{t+1},\ldots,u w_s\}+\{vw_{t+1},\ldots,v w_s\}.$ Note that $G-u\cong G-v.$ Hence $G_{s,t}\cong G_{t,s}.$ In view of Lemma~\ref{lem1.2}, one has $\lambda_\alpha(G_{t,s})>\lambda_\alpha(G_{s,t}),$ a contradiction to the fact that $\lambda_\alpha(G_{t,s})=\lambda_\alpha(G_{s,t})$. Hence $x_u\geqslant x_v.$ Applying Lemma~\ref{lem1.2} again, we have $\lambda_\alpha(G_{s,t})<\lambda_\alpha(G_{s+1,t-1}).$ This completes the proof.
\end{proof}
\begin{lem}\label{lem3.5}
If $s\geqslant t\geqslant 1$ and $s+t=n-2\geqslant 3$, then $\lambda_\alpha(G_1(s,t))>\lambda_\alpha(G_2(s-1,t))$ for $\alpha\in[0,1).$
\end{lem}
\begin{proof}
Note that $s\geqslant t\geqslant 1$ and $s+t=n-2\geqslant 3$. In what follows, we consider firstly that $t=1.$  Let ${\bf x}:=(x_1,x_2,\ldots,x_n)^T$ be a positive eigenvector of $A_\alpha(G_2(s-1,1))$ corresponding to $\lambda_\alpha(G_2(s-1,1))$.

If $s=2$, then $G_{2}(1,1)\cong P_5.$ Together with Lemma \ref{lem2.1}, our result follows obviously.

If $s=3$, then $G_2(2,1)\cong T_{1,1,3}.$ If $\alpha\in[0,s_1],$ then by Lemma \ref{lem1.3} one has $\lambda_\alpha(G_1(3,1))>2\geqslant \lambda_\alpha(G_2(2,1)).$ If $\alpha\in(s_1,1),$ then $\lambda_\alpha(G_2(2,1))>2.$ Based on Lemma \ref{lem3.3}, we have $x_{1}>x_{2}.$ In view of Lemma \ref{lem1.2}, we obtain $\lambda_\alpha(G_1(3,1))>\lambda_\alpha(G_2(2,1)).$

If $s\geqslant 4,$ then by Lemma~\ref{lem1.3} one has $\lambda_\alpha(G_2(s-1,1))>2$. In view of Lemma \ref{lem3.3}, we have $x_{1}>x_{2}.$ Together with Lemma \ref{lem1.2}, we obtain $\lambda_\alpha(G_1(s,1))>\lambda_\alpha(G_2(s-1,1)).$

Now we consider that $t\geqslant 2.$ Then let $G_1^*(s,t)$ be a graph obtained from $G_1(s,t)$ by subdividing the edge $v_1v_2.$ By Lemma~\ref{lem3.1}, one has $\lambda_\alpha(G_1^*(s,t))\leqslant \lambda_\alpha(G_1(s,t)).$ Notice that $G_2(s-1,t)$ is a proper subgraph of $G_1^*(s,t).$ Based on Lemma~\ref{lem3.2}, we have $\lambda_\alpha(G_2(s-1,t))<\lambda_\alpha(G_1^*(s,t))\leqslant\lambda(G_1(s,t)),$ as required.
\end{proof}

Notice that $G_1(s_1,t_1),\,G_2(s_2,t_2)\in \mathcal{G}_{n,n-2}.$ Now, we prove Theorem \ref{thm3.1}.
\begin{proof}[\bf Proof of Theorem \ref{thm3.1}]\
Let $G$ be a graph with  minimum $A_\alpha$-index in $\mathcal{G}_{n,n-2}.$ Assume, without loss of generality, that $I=\{v_3, v_4,\ldots, v_n\}$ is a maximum independent set of $G$. Notice that $\alpha(G)=n-2.$ Then each vertex of $I$ is adjacent to at least one of $\{v_1,v_2\}$ and $\min\{|N_G(v_1)\setminus \{v_2\}|,|N_G(v_2)\setminus \{v_1\}|\}\geqslant 1.$

If $(N_G(v_1)\setminus \{v_2\})\cap (N_G(v_2)\setminus \{v_1\})=\emptyset$, then by Lemma \ref{lem3.2} one has $G\cong G_1(s_1,t_1),$ where $s_1\geqslant t_1\geqslant1$ and $s_1+t_1=n-2.$ By Lemma \ref{lem3.5}, we have $\lambda_\alpha(G)>\lambda_\alpha(G_2(s_1-1, t_1)),$ which contradicts the minimality of $G.$

If $(N_G(v_1)\setminus \{v_2\})\cap (N_G(v_2)\setminus \{v_1\})\neq\emptyset$, then by Lemma \ref{lem3.2} one has $G\cong G_2(s_2,t_2)$ with {$\max\{s_2,t_2\}\geqslant 1$ and} $s_2+t_2=n-3.$ Together with Lemma \ref{lem3.4}, we have $G\cong G_2\left(\left\lceil\frac{n-3}{2}\right\rceil,\left\lfloor\frac{n-3}{2}\right\rfloor\right).$ This completes the proof.
\end{proof}

Next, we focus on proving Theorem \ref{thm3.2}. The following two lemmas will be needed in our proofs. It is straightforward to check that $H_1(s_1, t_1, k_1), H_2(s_2, t_2, k_2), H_3(s_3, t_3, k_3), H_4(s_4, t_4,k_4)\in \mathcal{G}_{n,n-3}.$
\begin{lem}\label{lem3.6}
Let $G$ be in $\mathcal{G}_{n,n-3}\,(n\geqslant 9)$ having the minimum $A_\alpha$-index. If $\alpha\in[0,1),$ then $G\cong H_3(s_3, t_3, k_3)$ for some {integers $s_3, t_3, k_3$ with $|s_3-k_3|\leqslant 1$ and} $s_3+t_3+k_3=n-5.$
\end{lem}
\begin{proof}
We firstly prove that $G\cong H_3(s_3, t_3, k_3)$ for some integers $s_3, t_3,k_3$ with $s_3+t_3+k_3=n-5.$
Assume, without loss of generality, that $S=\{v_4,v_5,\ldots,v_n\}$ is a maximum independent set of $G$. Hence each vertex of $S$ must adjacent to at least one of $\{v_1,v_2,v_3\}.$ Denote $N_G^*(v_1)=N_G(v_1)\setminus \{v_2,v_3\},$ $N_G^*(v_2)=N_G(v_2)\setminus \{v_1,v_3\},$ and $N_G^*(v_3)=N_G(v_3)\setminus \{v_1,v_2\}.$ Notice that $\alpha(G)=n-3.$ It is straightforward to check that $\min\left\{|N_G^*(v_1)|,|N_G^*(v_2)|,|N_G^*(v_3)|\right\}\geqslant 1.$ Now, we distinguish the following four cases.

{\bf Case 1.}\ $N_G^*(v_1)\cap N_G^*(v_2)= N_G^*(v_2)\cap N_G^*(v_3)= N_G^*(v_1)\cap N_G^*(v_3)=\emptyset.$  In this case, based on Lemma~\ref{lem3.2}, we have $G\cong H_1(s_1,t_1,k_1)$ with {$\min\{s_1, t_1, k_1\}\geqslant1$ and} $s_1+t_1+k_1=n-3.$ Without loss of generality, we assume that $s_1\geqslant k_1.$

If $s_1\geqslant 2,$ then let $H_1^*(s_1,t_1,k_1)$ be the graph obtained from $H_1(s_1,t_1,k_1)$ by subdividing the edge $v_1v_2.$ By Lemma~\ref{lem3.1}, we have $\lambda_\alpha(H_1^*(s_1,t_1,k_1))\leqslant \lambda_\alpha(H_1(s_1,t_1,k_1)).$ Note that $H_2(s_1-1,t_1,k_1)$ is a proper subgraph of $H_1^*(s_1,t_1,k_1).$ By Lemma \ref{lem3.2}, one has $\lambda_\alpha(H_1(s_1,t_1,k_1))\geqslant \lambda_\alpha(H_1^*(s_1,t_1,k_1))>\lambda_\alpha(H_2(s_1-1,t_1,k_1)),$ which contradicts the minimality of $G$.

If $s_1=k_1=1,$ then $t_1\geqslant 4.$ In view of Lemma \ref{lem1.3}, we have $\lambda_\alpha(H_1(s_1,t_1,k_1))>2.$ Let ${\bf x}:=(x_1,x_2,\ldots,x_n)^T$ be a positive eigenvector of $A_\alpha(H_2(s_1,t_1-1,k_1))$ corresponding to $\lambda_\alpha(H_2(s_1,t_1-1,k_1))$. By Lemma~\ref{lem3.3}, one has $x_{2}>x_{1}.$ Together with Lemma \ref{lem1.2}, we obtain $\lambda_\alpha(H_1(s_1,t_1,k_1))>\lambda_\alpha(H_2(s_1,t_1-1,k_1)),$ a contradiction.

{\bf Case 2.}\ $N_G^*(v_1)\cap N_G^*(v_2)\neq \emptyset,$ and $N_G^*(v_2)\cap N_G^*(v_3)=N_G^*(v_1)\cap N_G^*(v_3)=\emptyset.$ In this case, in view of Lemma~\ref{lem3.2}, we have $G\cong H_2(s_2,t_2,k_2)$ with {$\max\{s_2,t_2\}\geqslant 1,$} $k_2\geqslant1$ and $s_2+t_2+k_2=n-4.$ {If $k=1,$ then $G\cong H_3(s_3,t_3,0),$ as desired. If $k\geqslant 2,$ then by} a  similar discussion as that of Case 1, we can obtain a contradiction.

{\bf Case 3.}\ $N_G^*(v_1)\cap N_G^*(v_2)\neq \emptyset,$ $N_G^*(v_2)\cap N_G^*(v_3)\neq \emptyset$ and $N_G^*(v_1)\cap N_G^*(v_3)= \emptyset.$ In this case, by Lemma \ref{lem3.2}, we have $G\cong H_3(s_3,t_3,k_3)$ with {$\max\{s_3,t_3,k_3\}\geqslant 1$ and} $s_3+t_3+k_3=n-5.$

{\bf Case 4.}\ $N_G^*(v_1)\cap N_G^*(v_2)\neq \emptyset,$ $N_G^*(v_2)\cap N_G^*(v_3)\neq \emptyset$ and $N_G^*(v_1)\cap N_G^*(v_3)\neq \emptyset.$  If $N_G^*(v_1)\cap N_G^*(v_2)\cap N_G^*(v_3)= \emptyset,$ then there exists a vertex $w\in N_G^*(v_1)\cap N_G^*(v_2)$ and $w\not\in N_G^*(v_3).$ Without loss of generality, assume that $|N_G^*(v_1)|\geqslant |N_G^*(v_2)|\geqslant |N_G^*(v_3)|.$ Clearly, $|N_G^*(v_1)|\geqslant2.$ Let $G'=G-v_1w.$ It is routine to check that $G'\in \mathcal{G}_{n,n-3}.$ It follows from Lemma \ref{lem3.2} that $\lambda_\alpha(G')<\lambda_\alpha(G),$ a contradiction. Therefore, $N_G^*(v_1)\cap N_G^*(v_2)\cap N_G^*(v_3)\neq \emptyset.$

In this case, in view of Lemma \ref{lem3.2}, we have $G\cong H_4(s_4,t_4,k_4)$ with $\min\{s_4,\,t_4\}\geqslant1$ and $s_4+t_4+k_4=n-4.$ Let ${\bf x}:=(x_1,x_2,\ldots,x_n)^T$ be a positive eigenvector of $A_\alpha (H_3(s_4,t_4-1,k_4))$ corresponding to $\lambda_\alpha (H_3(s_4,t_4-1,k_4)).$ If $x_4\geqslant x_5,$ then by Lemma \ref{lem1.2}, one has $\lambda_\alpha(H_3(s_4,t_4-1,k_4))<\lambda_\alpha(H_3(s_4,t_4-1,k_4)-v_3v_5+v_3v_4).$ Notice that $H_3(s_4,t_4-1,k_4)-v_3v_5+v_3v_4\cong H_4(s_4,t_4,k_4).$ Hence $\lambda_\alpha(H_3(s_4,t_4-1,k_4))<\lambda_\alpha(G),$ a contradiction. If $x_4<x_5,$ then by Lemma \ref{lem1.2}, one has $\lambda_\alpha(H_3(s_4,t_4-1,k_4))<\lambda_\alpha(H_3(s_4,t_4-1,k_4)-v_1v_4+v_1v_5)=\lambda_\alpha(G),$ a contradiction.

Combining with Cases 1-4, we obtain $G\cong H_3(s_3, t_3, k_3)$ for some integers $s_3, t_3,k_3$ with $s_3+t_3+k_3=n-5.$ Together with Lemma \ref{lem3.4}, our result follows immediately.
\end{proof}
Notice that if $H\cong H_3(s_3,t_3,k_3)$ for some integers $s_3,t_3,k_3$ with $\min\{s_3,t_3,k_3\}\geqslant 1$ and $s_3+t_3+k_3=n-5,$ then
$$
\pi_H=(N_H(v_1)\setminus\{v_4\})\cup (N_H(v_3)\setminus\{v_5\})\cup \{v_1\}\cup \{v_3\}\cup \{v_4\}\cup \{v_5\}\cup \{v_2\}\cup
(N_H(v_2)\setminus \{v_4,v_5\})
$$
is an equitable partition of $V_H.$ Denote by $Q_{\pi_H}$ the quotient matrix of $A_\alpha(H)$ corresponding to the partition $\pi_H.$
\begin{lem}\label{lem3.11}
Let $G$ be a graph among $\mathcal{G}_{n,n-3}$ with $n\geqslant 9$ having the minimum $A_\alpha$-index. If $\alpha\in[\frac{1}{2},1)$ and $n\not\equiv2\pmod{3},$ or $\alpha\in[\frac{11}{20},1)$ and $n\equiv2\pmod{3},$ then $G\cong H_3(s, t, s)$ for some $s, t$ with $2s+t=n-5.$
\end{lem}
\begin{proof}
By Lemma \ref{lem3.6}, it suffices to prove that $G\not\cong H_3(s,t,s-1)$ {for some nonnegative integers $s,t$ with $s\geqslant 1$ and $2s+t=n-4.$} Suppose to the contrary that $G\cong H_3(s,t,s-1).$ We proceed by considering the following cases.

{\bf Case 1.} $s\leqslant \frac{n-3}{3}.$ Then $s\leqslant t+1.$ {Note that $t\geqslant 1.$ If $s=t+1,$ then $G\cong H_3(t+1,t,t).$ In view of Lemma \ref{lem3.7}, one has $\lambda_\alpha(H_3(t+1,t,t))>\alpha(t+2).$ By \textit{Mathematica} \cite{math}, we obtain that if $x> \alpha(t+2),$ then
$$
  \det(xI-Q_{\pi_{H_3(t+1,t-1,t+1)}})-\det(xI-Q_{\pi_{H_3(t+1,t,t)}})=(1-\alpha)^2(1+(x-2)\alpha)^2(x-2\alpha)(x-(t+2)\alpha)>0.
$$
Together with Lemmas \ref{lem3.9} and \ref{lem4.01}, we know that the $A_\alpha$-index of $H_3(t+1,t-1,t+1)$ is less than that of $H_3(t+1,t,t),$ a contradiction.}

Next, we consider that $s\leqslant t.$ By Lemmas \ref{lem3.7} and \ref{lem3.8}, one has
$$
  \lambda_\alpha(H_3(s, t, s-1))>\alpha (t+2)+\frac{(1-\alpha)^2}{\alpha},\ \ \lambda_\alpha(H_3(s,t-1,s))\leqslant \alpha (t+1)+2(1-\alpha).
$$
Notice that
$$
  \left(\alpha (t+2)+\frac{(1-\alpha)^2}{\alpha}\right)- \left(\alpha (t+1)+2(1-\alpha)\right)=\frac{\alpha^2+(1-\alpha)^2-2\alpha(1-\alpha)}{\alpha}=\frac{(2\alpha-1)^2}{\alpha}\geqslant 0.
$$
Hence $\lambda_\alpha(H_3(s,t-1,s))<\lambda_\alpha(H_3(s, t, s-1)),$ a contradiction.

{\bf Case 2.} $s\geqslant \frac{n-1}{3}.$ Then $s\geqslant t+3.$ By Lemmas \ref{lem3.7} and \ref{lem3.8}, we have
$$
  \lambda_\alpha(H_3(s, t, s-1))>\alpha (s+1)+\frac{(1-\alpha)^2}{\alpha},\ \ \lambda_\alpha(H_3(s-1,t+1,s-1))\leqslant \alpha s+2(1-\alpha).
$$
Note that
$$
  \left(\alpha (s+1)+\frac{(1-\alpha)^2}{\alpha}\right)- \left(\alpha s+2(1-\alpha)\right)=\frac{\alpha^2+(1-\alpha)^2-2\alpha(1-\alpha)}{\alpha}=\frac{(2\alpha-1)^2}{\alpha}\geqslant 0.
$$
Hence $\lambda_\alpha(H_3(s-1,t+1,s-1))<\lambda_\alpha(H_3(s, t, s-1)),$ a contradiction.

{\bf Case 3.} $\frac{n-3}{3}<s< \frac{n-1}{3}.$ If $n=3m$ or $n=3m+1,$ then there is no integer satisfying $\frac{n-3}{3}<s< \frac{n-1}{3}.$ If $n=3m+2,$ then $s=m$ and $G\cong H_3(m,m-2,m-1).$ Notice that $m\geqslant 3$ and $\alpha\in [\frac{11}{20},1)$ in this case. For convenience, we put $\tilde{H}:=H_3(m,m-2,m-1)$ and $\hat{H}:=H_3(m,m-3,m).$ In what follows, we show that $\lambda_\alpha(\tilde{H})>\lambda_\alpha(\hat{H}).$ Further on we need the following claim. %
\begin{claim}\label{c1}
  If $n=3m+2$ and $\alpha\in[\frac{11}{20},1),$ then $\alpha(m+1){<} \lambda_\alpha(\hat{H})< \alpha m+1.$
\end{claim}
\begin{proof}[\bf Proof of Claim \ref{c1}.]\
In view of Lemma \ref{lem3.7}, one has $\lambda_\alpha(\hat{H}){>} \alpha(m+1).$ Notice that
$$
{\hat{\pi}}=((N_{\hat{H}}(v_1)\cup N_{\hat{H}}(v_3))\setminus\{v_4,v_5\})\cup\{v_1,v_3\}\cup\{v_4,v_5\}\cup \{v_2\}\cup (N_{\hat{H}}(v_2)\setminus \{v_4,v_5\})
$$
is an equitable partition of $V_{\hat{H}}.$ It is straightforward to check that the quotient matrix of $A_{\alpha}(\hat{H})$ corresponding to the partition $\hat{\pi}$ is
\begin{equation*}
    {A_\alpha(\hat{H})_{\hat{\pi}}}=\left(
  \begin{array}{ccccc}
    \alpha & 1-\alpha & 0 & 0 & 0 \\
    (1-\alpha)m & \alpha(m+1) & 1-\alpha & 0 & 0 \\
    0 & 1-\alpha & 2\alpha & 1-\alpha & 0 \\
    0 & 0 & 2(1-\alpha) & \alpha(m-1) & (1-\alpha)(m-3) \\
    0 & 0 & 0 & 1-\alpha & \alpha \\
  \end{array}
\right).
\end{equation*}
By a direct calculation, we have
\begin{align*}
    \det\left(xI-A_\alpha(\hat{H})_{\hat{\pi}}\right)=&x^5-2\alpha(2+m)x^4+(4\alpha^2-2m +4\alpha m+6\alpha^2m+\alpha^2m^2)x^3+(-5 \alpha+10 \alpha^2-3\alpha^3\\
     &+6\alpha m-12\alpha^2m-4\alpha^3 m+2\alpha m^2-4\alpha^2m^2-2\alpha^3 m^2)x^2+(-3+12\alpha-8\alpha^2-8\alpha^3\\
     &+2\alpha^4-\alpha^2m+2\alpha^3 m+3\alpha^4 m+m^2-4\alpha m^2+8\alpha^3 m^2)x+(-2\alpha+8\alpha^2-8\alpha^3)m^2\\
     &+ (3\alpha-12\alpha^2+ 15 \alpha^3-6 \alpha^4)m-2\alpha^4+13\alpha^3-12\alpha^2+3 \alpha
\end{align*}
and
\begin{align*}
 \det\left(\left(\alpha m+1\right)I-A_\alpha(\hat{H})_{\hat{\pi}}\right)=&(1-\alpha)\big(-2m\alpha^4+(6 m^3-16m^2+15m)\alpha^3-(5m^3-22m^2+22m+2)\alpha^2\\
 &+(m^3-9m^2+11m+4)\alpha+m^2-2(m+1)\big),
\end{align*}
here $I$ denotes the identity matrix. Let
$$
l(x)=\det(xI-A_\alpha(\hat{H})_{\hat{\pi}})
$$
be a real function {of} $x$ for $x\in [\alpha m+1,+\infty),$ and let
$$
f(\alpha)=-2m{\alpha}^4+(6 m^3-16m^2+15m){\alpha}^3-(5m^3-22m^2+22m+2){\alpha}^2+(m^3-9m^2+11m+4){\alpha}+m^2-2(m+1)
$$
be a real function {of} ${\alpha}$ for ${\alpha}\in [\frac{1}{2},1).$

Firstly, we show that $\det\left(\left(\alpha m+1\right)I-A_\alpha(\hat{H})_{\hat{\pi}}\right)>0$ if $\alpha\in[\frac{11}{20},1).$ It is routine to check that
\begin{align*}
    &f^{(1)}({\alpha})=-8m{\alpha}^3+3(6 m^3 -16 m^2+15m){\alpha}^2-2(5 m^3 -22 m^2+22 m+2){\alpha}+ m^3-9 m^2+11 m +4,\\
    &f^{(2)}({\alpha})=-24m{\alpha}^2+6(6 m^3 -16 m^2+15m){\alpha}-2(5 m^3 -22 m^2+22 m+2),\\
    &f^{(3)}({\alpha})=-48m{\alpha}+6(6 m^3 -16 m^2+15m),\ \ \ \text{and}\ \ \ f^{(4)}({\alpha})=-48m<0,
\end{align*}
where $f^{(i)}({\alpha})$ denotes the $i$-th derivative function of $f({\alpha}).$ Hence $f^{(3)}({\alpha})$ is decreasing in the interval $[\frac{1}{2},1).$ Notice that $m\geqslant 3.$ Therefore,
$$
  f^{(3)}({\alpha})\geqslant f^{(3)}(1)=6m(6m^2-16m+7)>0,
$$
That is, $f^{(2)}({\alpha})$ is increasing in the interval $[\frac{1}{2},1).$ Thus,
$$
  f^{(2)}({\alpha})\geqslant f^{(2)}(\frac{1}{2})=m(8m^2-4m-5)-4>0. $$
Then $f^{(1)}({\alpha})$ is increasing in the interval $[\frac{1}{2},1).$ Hence,
$$
  f^{(1)}({\alpha})\geqslant f^{(1)}(\frac{1}{2})=\frac{1}{4}(2m^3+4m^2-3m+8)>0.
$$
That is, $f({\alpha})$ is increasing in the interval $[\frac{1}{2},1).$ Hence, if $\alpha\in[\frac{11}{20},1)$ and $m\geqslant 3,$ then
$$
  f({\alpha})\geqslant f(\frac{11}{20})={\frac{2860 m^{3} + 3440m^{2} -23391m -32400}{80000}\geqslant \frac{5607}{80000}>0}.
$$
It follows that for $\alpha\in[\frac{11}{20},1),$
\[\label{eq:3.2}
\det\left(\left(\alpha m+1\right)I-A_\alpha(\hat{H})_{\hat{\pi}}\right)>0.
\]

{Now, we show that $l(x)$ is an increasing function {of $x$} for $x\in[\alpha m+1,+\infty).$ It is straightforward to check that
\begin{align*}
    l^{(1)}(x)=&5x^4-8\alpha(2+m)x^3+3(4\alpha^2-2m +4\alpha m+6\alpha^2m+\alpha^2m^2)x^2+2(-5 \alpha+10 \alpha^2-3\alpha^3\\
     &+6\alpha m-12\alpha^2m-4\alpha^3 m+2\alpha m^2-4\alpha^2m^2-2\alpha^3 m^2)x+(-3+12\alpha-8\alpha^2-8\alpha^3\\
     &+2\alpha^4-\alpha^2m+2\alpha^3 m+3\alpha^4 m+m^2-4\alpha m^2+8\alpha^3 m^2),\\
    l^{(2)}(x)=&20x^3-24\alpha(2+m)x^2+6(4\alpha^2-2m +4\alpha m+6\alpha^2m+\alpha^2m^2)x+2(-5 \alpha+10 \alpha^2-3\alpha^3\\
     &+6\alpha m-12\alpha^2m-4\alpha^3 m+2\alpha m^2-4\alpha^2m^2-2\alpha^3 m^2),\\
     l^{(3)}(x)=&60x^2-48\alpha(2+m)x+6(4\alpha^2-2m +4\alpha m+6\alpha^2m+\alpha^2m^2), \ \text{and}\\
     l^{(4)}(x)=&120x-48\alpha(2+m)\geqslant 120(\alpha m+1)-48\alpha(2+m)>0.
\end{align*}
Hence, $l^{(3)}(x)$ is an increasing function {of $x$} for $x\in[\alpha m+1,+\infty).$ Thus,
$$
  l^{(3)}(x)\geqslant l^{(3)}(\alpha m+1)={6( \alpha^{2}(3 m^{2} -10m + 4) +16(m - 1)\alpha -2m + 10)>0}.
$$
So, $l^{(2)}(x)$ is an increasing function {of $x$} for $x\in[\alpha m+1,+\infty).$ Thus,
\begin{align*}
  l^{(2)}(x)\geqslant &l^{(2)}(\alpha m+1)={2(\alpha^{3}(m^{3} -8m^{2} + 8m-3) + \alpha^{2}(17m^{2}- 42m +22)}\\
  &{-\alpha(4m^{2}-36m+29) - 6m + 10)>0.}
\end{align*}
Hence $l^{(1)}(x)$ is an increasing function {of $x$} for $x\in[\alpha m+1,+\infty).$ Thus,
\begin{align*}
  l^{(1)}(x)\geqslant &l^{(1)}(\alpha m+1)={\alpha^4 (- 2 m^3 + 4 m^2- 3 m+2)+2\alpha^3 (3 m^3- 16 m^2 +19 m -7)}\\
  &{+ \alpha^2 (-2m^3+ 37 m^2- 65 m+24)- 2\alpha( 6 m^2 - 18 m +7)+m(m-6)+2>0.}
\end{align*}
That is to say, $l(x)$ is an increasing function {of $x$} for $x\in[\alpha m+1,+\infty).$

Together with \eqref{eq:3.2}, we know that the largest eigenvalue of $A_\alpha(\hat{H})_{\hat{\pi}}$ is less than $\alpha m+1.$ By Lemma~\ref{lem4.01}, one obtains that $\lambda_\alpha(\hat{H})$ equals to the largest eigenvalue of $A_\alpha(\hat{H})_{\hat{\pi}}.$ Hence $\lambda_\alpha(\hat{H})< \alpha m+1$ if ${\alpha}\in [\frac{11}{20},1).$}
This completes the proof of Claim \ref{c1}.
\end{proof}

By Claim \ref{c1}, we know that $2<\alpha(m+1){<}\lambda_\alpha(\hat{H})< \alpha m+1.$ By \textit{Mathematica} \cite{math}, we obtain that for $\alpha\in[\frac{11}{20},1),$ $m\geqslant 3$ and $x\in {(}\alpha(m+1),\alpha m+1),$
\begin{align*}
  \det(xI-Q_{\pi_{\tilde{H}}})-\det(xI-Q_{\pi_{\hat{H}}})=&(1 + \alpha (x-2))^2 (x-2\alpha)\big(x^3 +\alpha (2 + 3 m + 2 (2 + m)x- (4 + m)x^2)\\
  &- (2 + m)x+\alpha^2 (2 m (x-3) + 3 x-4)+\alpha^3 m\big)<0.
\end{align*}

Together with Lemma \ref{lem3.9}, we derive that $\lambda_\alpha(\tilde{H})>\lambda_\alpha(\hat{H}),$ which contradicts the minimality of $\tilde{H}$.

Combining with Cases 1-3, we know that $G\cong H_3(s, t, s)$ for some $s, t$ with $2s+t=n-5.$ This completes the proof.
\end{proof}
Now, we are ready to prove Theorem \ref{thm3.2}.
\begin{proof}[\bf Proof of Theorem \ref{thm3.2}.]\
The first assertion of Theorem \ref{thm3.2} can be derived from Lemma \ref{lem3.6}. In what follows, we focus on proving the second assertion of this  theorem. That is, we consider the case that $\alpha\in[\frac{1}{2},1)$ if $n\not\equiv2\pmod{3}$ and $\alpha\in[\frac{11}{20},1)$ if $n\equiv2\pmod{3}.$

Choose $G$ in $\mathcal{G}_{n,n-3}$ such that its $A_\alpha$-index is as small as possible. According to Lemma~\ref{lem3.11}, one has
$G\cong H_3(s, t, s)$ for some $s, t$ with $2s+t=n-5.$ We proceed by considering the following two cases.

{\bf Case 1.}\ $s\geqslant t+1$. In fact, $s\leqslant t+3.$ Otherwise, $s\geqslant t+4.$ By Lemmas \ref{lem3.7} and \ref{lem3.8}, we have
$$
\lambda_\alpha(H_{3}(s,t,s))> \alpha(s+1)+\frac{(1-\alpha)^2}{\alpha}\ \text{and}\ \lambda_\alpha(H_{3}(s-1,t+2,s-1))\leqslant \alpha s+2(1-\alpha).
$$
Hence $\lambda_\alpha(H_{3}(s,t,s))>\lambda_\alpha(H_{3}(s-1,t+2,s-1)),$ a contradiction.  Therefore, 
\begin{equation*}
    G=
    \left\{
    \begin{aligned}
        &H_3(m-1,m-3,m-1),\ \ &\textrm{if $n=3m;$}\\
        &H_3(m-1,m-2,m-1),\ \ &\textrm{if $n=3m+1;$}\\
        &H_3(m,m-3,m),\ \ &\textrm{if $n=3m+2.$}
    \end{aligned}
    \right.
\end{equation*}




{\bf Case 2.}\ $s\leqslant t+1$. In fact, $s\geqslant t-1.$ Otherwise, $s\leqslant t-2.$ By Lemmas \ref{lem3.7} and \ref{lem3.8}, we have
$$
\lambda_\alpha(H_{3}(s,t,s))> \alpha(t+2)+\frac{(1-\alpha)^2}{\alpha}\ \text{and}\ \lambda_\alpha(H_{3}(s+1,t-2,s+1))\leqslant \alpha t+2(1-\alpha).
$$
Hence $\lambda_\alpha(H_{3}(s,t,s))>\lambda_\alpha(H_{3}(s+1,t-2,s+1)),$ a contradiction.  Therefore, 
\begin{equation*}
    G=
    \left\{
    \begin{aligned}
        &H_3(m-2,m-1,m-2),&\ \ \textrm{if $n=3m;$}\\
        &H_3(m-1,m-2,m-1),&\ \ \textrm{if $n=3m+1;$}\\
        &H_3(m-1,m-1,m-1),\ \ &\textrm{if $n=3m+2.$}
    \end{aligned}
    \right.
\end{equation*}




Now, we compare the {spectral radii of the extremal graphs} in Case 1 with those in Case 2.

If $n=3m,$ then by Lemmas \ref{lem3.7} and \ref{lem3.8}, one has $\lambda_\alpha(H_3(m-2, m-1, m-2))>\alpha(m+1)+\frac{(1-\alpha)^2}{\alpha}$ and $\lambda_\alpha(H_3(m-1, m-3, m-1))\leqslant \alpha m+2(1-\alpha).$ By a direct calculation, one has $\lambda_\alpha(H_3(m-2, m-1, m-2))>\lambda_\alpha(H_3(m-1, m-3, m-1)).$ Hence $G\cong H_3(m-1, m-3, m-1).$

If $n=3m+1,$ then $G\cong H_3(m-1,m-2,m-1).$

If $n=3m+2,$ then denote $H':=H_3(m-1,m-1,m-1)$ and $\hat{H}:=H_3(m,m-3,m).$ By a direct calculation, we have
\begin{align}\notag
    \det(xI- A_\alpha(H')_{\pi'})-\det(xI- A_\alpha(\hat{H})_{\hat{\pi}})=&(1+\alpha(x-2))(x-2\alpha)(1+m-2\alpha(1+m)\\
    &+\alpha(2+m)x-x^2),\label{eq:3.1}
\end{align}
where
$$
  \pi'=((N_{H'}(v_1)\cup N_{H'}(v_3))\setminus\{v_4,v_5\})\cup \{v_1,v_3\}\cup \{v_4,v_5\}\cup \{v_2\}\cup (N_{H'}(v_2)\setminus \{v_4,v_5\})
$$
and
$$
\hat{\pi}=((N_{\hat{H}}(v_1)\cup N_{\hat{H}}(v_3))\setminus\{v_4,v_5\})\cup\{v_1,v_3\}\cup\{v_4,v_5\}\cup \{v_2\}\cup (N_{\hat{H}}(v_2)\setminus \{v_4,v_5\})
$$
are equitable partitions of $H'$ and $\hat{H},$ respectively.

It is routine to check that $x_1:=\frac{1}{2}(\alpha(2 + m)+r)$ is the largest root of $1+m-2\alpha(1 + m)+\alpha(2+m)x-x^2=0,$ where $r=\sqrt{\alpha^2 (2 + m)^2 + 4 (1-2\alpha)(1+m)}.$ By \textit{Mathematica} \cite{math}, we obtain  that for $\alpha\in[\frac{11}{20},1)$ and $m\geqslant 3,$
$$
  \det(x_1I- A_\alpha(\hat{H})_{\hat{\pi}})=-(\alpha-1)^2(r+\alpha(4+m-2r+2\alpha(\alpha(2+m)+r-m-4)))<0.
$$
Therefore, $\lambda_\alpha(\hat{H})>x_1.$ In view of Lemma \ref{lem1.3}, one has $\lambda_\alpha(\hat{H})>2.$ Hence \eqref{eq:3.1} is less than $0$ if $x=\lambda_\alpha(\hat{H})$. Together with Lemma~\ref{lem3.9}, we have $\lambda_\alpha(H')>\lambda_\alpha(\hat{H}).$ Hence $G\cong H_3(m,m-3,m).$ This completes the proof.
\end{proof}
\section{\normalsize Proof of Theorem \ref{thm5.1}}
In this section, we prove Theorem \ref{thm5.1}, which characterizes the graph with  minimum $A_\alpha$-index among $\mathcal{G}_{n,2}$.  It is routine to check that if $G$ is a graph among $\mathcal{G}_{n,2}$, then $G^c$ is $K_3$-free. Choose $v_1\in V_G$ such that $d_{G^c}(v_1)=\Delta(G^c).$ Assume, without loss of generality, that $N_{G^c}(v_1)=\{v_{n-\Delta(G^c)+1},\ldots,v_{n-1},v_n\}.$ Hence $G^c[N_{G^c}(v_1)]$ is empty.

We prove the following structural lemma at first.
\begin{lem}\label{lem4.2}
Let $G$ be in $\mathcal{G}_{n,2}$  with  minimum $A_\alpha$-index and $n\geqslant 11$. If $\alpha\in [0,\frac{3}{4}{]},$ then $G^c$ is a bipartite graph.
\end{lem}
\begin{proof}
Notice that $G^c$ is $K_3$-free. Then the induced subgraph $G^c[N_{G^c}(v_1)]$ is an empty graph. In order to complete the proof, it is sufficient to show that $G^c[V_G\setminus N_{G^c}(v_1)]$ is empty. We proceed by distinguishing the parity on $n$.

{\bf Case 1.}\ $n=2t.$ {In this case, $t\geqslant 6.$} In view of Lemma \ref{lem4.1},  we have $|E_{G^c}|\leqslant t^2-1.$ That is to say, $|E_G|\geqslant t^2-t+1.$ We are to show that $|E_G|= t^2-t+1.$ In fact, if $|E_G|\geqslant t^2-t+2,$ then by Lemma~\ref{lem3.8} one has
$$
  \lambda_\alpha(G)\geqslant \frac{2|E_G|}{n}\geqslant t-1+\frac{2}{t},
$$
{and both of the equalities hold if and only if $G$ is regular and $|E_G|=t^2-t+2,$ i.e., $G$ is $(t-1+\frac{2}{t})$-regular, which is impossible since $t\geqslant 6.$ Hence $\lambda_\alpha(G)>t-1+\frac{2}{t}.$}

Notice that $\pi_1=(V_{F(t,t)}\setminus \{u,v\})\cup \{u,v\}$ is an equitable partition of $V_{F(t,t)}.$ It is routine to check that the quotient matrix of $A_\alpha(F(t,t))$ corresponding to the partition $\pi_1$ can be written as
\begin{equation*}
   {A_\alpha(F(t,t))_{\pi_1}}=\left(
       \begin{array}{cc}
         t-2+\alpha & 1-\alpha \\
         (1-\alpha)(t-1) & \alpha t+1-\alpha \\
       \end{array}
     \right).
\end{equation*}
Together with Lemma \ref{lem4.01}, we obtain
$$
\lambda\left(\det(xI- A_\alpha(F(t,t))_{\pi_1})\right)=\lambda_\alpha(F(t,t))=\frac{\alpha t+t-1+ \sqrt{(\alpha t+t-1)^2-4(\alpha t^2-\alpha t-1+\alpha)}}{2}.
$$
Now, we show that for $\alpha\in [0,\frac{3}{4}],$
$$
  t-1+\frac{2}{t}\geqslant \frac{\alpha t+t-1+ \sqrt{(\alpha t+t-1)^2-4(\alpha t^2-\alpha t-1+\alpha)}}{2}.
$$
It suffices to show
$$
  (t-1+\frac{4}{t}- \alpha t)^2\geqslant (\alpha t+t-1)^2-4(\alpha t^2-\alpha t-1+\alpha),
$$
which is equivalent to
$$
  (1-\alpha)t^2-2t+4\geqslant 0.
$$
Notice that $4-16(1-\alpha)^2\leqslant 0$ if $\alpha\in [0,\frac{3}{4}].$ Hence $(1-\alpha)t^2-2t+4\geqslant 0$ holds for $\alpha\in [0,\frac{3}{4}].$ Therefore, $\lambda_\alpha(G)>\lambda_\alpha(F(t,t)),$ a contradiction to the minimality of $G.$ Hence $|E_G|= t^2-t+1$ and $|E_{G^c}|=t^2-1.$

Suppose that $G^c[V_G\setminus N_{G^c}(v_1)]$ is not empty. Then there exists an edge, say $v_iv_j$, in $E_{G^c}$ with {$v_i,v_j\in V_G\setminus N_{G^c}(v_1).$} Recall that $G^c$ is $K_3$-free. Hence, the size of $\{xy\in E(G^c): x\in \{v_i,v_j\}, y\in N_{G^c}(v_1)\}$  is less than or equal to $\Delta(G^c)$. {Denote $W:=V_G\setminus(N_{G^c}(v_1)\cup\{v_i,v_j\}).$} Therefore,
\begin{align*}
    t^2-1&=|E_{G^c}|\leqslant\sum_{{v\in W}}d_{G^c}(v)+\Delta(G^c)+1\leqslant \Delta(G^c)(2t-\Delta(G^c)-2)+\Delta(G^c)+1\\
    &=-(\Delta(G^c)-t)^2+t^2-\Delta(G^c)+1\leqslant t^2-2,
\end{align*}
a contradiction.

{\bf Case 2.}\ $n=2t+1.$ {In this case, $t\geqslant 5.$} By Lemma \ref{lem4.1}, one has $|E_{G^c}|\leqslant t^2+t-1.$ Hence $|E_G|\geqslant t^2+1.$ We claim that $|E_G|\leqslant t^2+t-2$. In fact, if $|E_G|\geqslant t^2+t-1,$ then by Lemma~\ref{lem3.8} we have $\lambda_\alpha(G)\geqslant \frac{2|E_G|}{n}\geqslant t+\frac{1}{2}-\frac{5}{4t+2}.$ Note that $\pi_2=(V_{K_t}\setminus \{u\})\cup \{u\}\cup \{v\}\cup (V_{K_{t+1}}\setminus \{v\})$ is an equitable partition of $V_{F(t,t+1)}.$ It is straightforward to check that the quotient matrix of $A_\alpha(F(t,t+1))$ corresponding to the partition $\pi_2$ can be written as
\begin{equation*}
    {A_\alpha(F(t,t+1))_{\pi_2}}=\left(
                      \begin{array}{cccc}
                        t-2+\alpha & 1-\alpha & 0 & 0 \\
                        (1-\alpha)(t-1) & \alpha t & 1-\alpha & 0 \\
                        0 & 1-\alpha & \alpha(t+1) & (1-\alpha)t \\
                        0 & 0 & 1-\alpha & t-1+\alpha \\
                      \end{array}
                    \right).
\end{equation*}
We proceed by considering the following claim, whose detailed proof is given in the Appendix.
\begin{claim}\label{c3}
If $\alpha\in [0,\frac{3}{4}],$ then $\lambda_\alpha(F(t,t+1))<t+\frac{1}{4}.$
\end{claim}

In view of Claim \ref{c3}, we have $\lambda_\alpha(F(t,t+1))<t+\frac{1}{4}< t+\frac{1}{2}-\frac{5}{4t+2}\leqslant \lambda_\alpha(G)$ if $t\geqslant 5,$ which contradicts the minimality of $G$. Therefore, $|E_G|\leqslant t^2+t-2$, and so $|E_{G^c}|\geqslant t^2+2.$

Suppose that $G^c[V_G\setminus N_{G^c}(v_1)]$ is not empty. Then there exists an edge, say $v_rv_l$, in $E_{G^c}$ with {$v_r,v_l\in V_G\setminus N_{G^c}(v_1).$}  Recall that $G^c$ is $K_3$-free. Hence, the size of $\{xy\in E(G^c): x\in \{v_r,v_l\}, y\in N_{G^c}(v_1)\}$  is no more than  $\Delta(G^c)$. {Denote $W':=V_G\setminus(N_{G^c}(v_1)\cup\{v_r,v_l\}).$} Therefore,
\begin{align*}
    t^2+2&{\leqslant} |E_{G^c}|\leqslant \sum_{{v\in W'}}d_{G^c}(v)+\Delta(G^c)+1\leqslant \Delta(G^c)(2t+1-\Delta(G^c)-2)+\Delta(G^c)+1\\
    &=-(\Delta(G^c)-t)^2+t^2+1\leqslant t^2+1,
\end{align*}
a contradiction.

By Cases 1 and 2, the result follows immediately.
\end{proof}

\begin{proof}[\bf Proof of Theorem \ref{thm5.1}]
By Lemmas \ref{lem3.2} and \ref{lem4.2}, one has $G\cong F(s,t)$ for some $s\geqslant t$ and $s+t=n.$ If $s=t$ or $t+1,$ then our result holds. So, in what follows, we assume that $s\geqslant t+2.$ Note that $\pi=(V_{K_s}\setminus \{u\})\cup \{u\}\cup \{v\}\cup (V_{K_t}\setminus \{v\})$ is an equitable partition of $V_{F(s,t)}.$ So the quotient matrix of $A_\alpha(F(s,t))$ corresponding to the partition $\pi$ is
\begin{equation*}
    {A_\alpha(F(s,t))_{\pi}}=\left(
                      \begin{array}{cccc}
                        s-2+\alpha & 1-\alpha & 0 & 0 \\
                        (1-\alpha)(s-1) & \alpha s & 1-\alpha & 0 \\
                        0 & 1-\alpha & \alpha t & (1-\alpha)(t-1) \\
                        0 & 0 & 1-\alpha & t-2+\alpha \\
                      \end{array}
                    \right).
\end{equation*}
Now, we distinguish the following two cases.

{\bf Case 1.}\ $s= t+2.$ By a direct calculation, we have
\[\label{eq:4.4}
  \det((t+1)I- A_\alpha(F(t,t+2))_{\pi})=(1-\alpha)\big(\alpha(3+2\alpha(t-1)-2t)-3\big)<0.
\]
The detailed proof for the last inequality in \eqref{eq:4.4} is given in the Appendix.

Together with Lemma \ref{lem4.01}, we obtain  $\lambda_\alpha(F(t+2,t))>t+1.$ Based on Lemma \ref{lem3.8}, one has
$$
\lambda_\alpha(F(t+1,t+1))\leqslant t+1<\lambda_\alpha(F(t+2,t)),
$$
a contradiction.

{\bf Case 2.}\ $s\geqslant t+3$. In this case, we consider firstly that $\alpha\in [0,\frac{1}{2}].$
By a direct calculation, one has
\begin{align}\notag
    \det((s-1-\alpha)I- A_\alpha(F(s,t))_{\pi})=&t-1-s+\alpha\big(\alpha^3(5+3t+s(3+t))+\alpha^2(s(t^2-2t-11)\\\notag
    &+3(t^2+t+2)-s^2(t+3))+\alpha(s^3+6s^2-s(t^2+4t+1)\\\label{eq:4.2}
    &-(t-3)t-12)+s(1+s)(t-s)+4s-4t+6\big)<0.
\end{align}
The detailed proof for the last inequality in \eqref{eq:4.2} is given in the Appendix.

Together with Lemma \ref{lem4.01}, we obtain  $\lambda_\alpha(F(s,t))>s-1-\alpha.$ In view of Lemma \ref{lem3.8}, one has
$$
  \lambda_\alpha(F(s-1,t+1))\leqslant s-1-\alpha<\lambda_\alpha(F(s,t)),
$$
a contradiction.

Now we consider that $\alpha\in [\frac{1}{2},\frac{3}{4}].$
By a direct calculation, we have
\[\label{eq:4.3}
  \det((s-2+\alpha)I- A_\alpha(F(s,t))_{\pi})=(\alpha-1)^2 (s-1)\big((s-1)(1-s+t)+\alpha^2(t-1)+\alpha(s-t-2)(t-1)\big)<0.
\]
The detailed proof for the last inequality in \eqref{eq:4.3} is given in the Appendix.

Together with Lemma \ref{lem4.01}, we obtain $\lambda_\alpha(F(s,t))>s-2+\alpha.$ In view of Lemma \ref{lem3.8}, one has
$$
\lambda_\alpha(F(s-1,t+1))\leqslant s-2+\alpha<\lambda_\alpha(F(s,t)),
$$
a contradiction.

Together with Cases 1 and 2, one has $G\cong F(\left\lceil\frac{n}{2}\right\rceil,\left\lfloor\frac{n}{2}\right\rfloor).$ This completes the proof.
\end{proof}
\section{\normalsize Proofs of Theorems \ref{thm1.5} and \ref{thm1.6}}\setcounter{equation}{0}
In this section, we determine the unique $n$-vertex tree (resp. bipartite graph and simple graph) with given independence number having the maximum $A_\alpha$-index.

Notice that if $G$ is a connected bipartite graph of order $n$, then the independence number $i(G)$ and the matching number $\mu(G)$ satisfy
\[\label{eq:3.01}
\text{$i(G) \geqslant \left\lceil\frac{n}{2}\right\rceil$\ \  and\ \ $i(G)+\mu(G)=n$.}
\]

\begin{proof}[\bf Proof of Theorem \ref{thm1.5}]
Choose an {$n$-vertex bipartite graph $G$ with independence number $i$} \ ($\left\lceil\frac{n}{2}\right\rceil\leqslant i\leqslant n-1$) such that $\lambda_\alpha(G)$ is as large as possible. Assume that $G$ has color classes $U,$ $V$ and $|U|\geqslant |V|.$ Hence $i \geqslant |U|.$

If $i =|U|,$ then by Lemma \ref{lem3.2} one has $G\cong K_{i,n-i},$ as desired. If $i>|U|$, then $G\not\cong K_{|U|,|V|}.$ Notice that $\alpha\in [\frac{1}{2},1).$ Applying Lemma \ref{lem3.2} again, we have
\begin{align*}
\lambda_\alpha(G)&<\lambda_\alpha(K_{|U|,|V|})=\frac{\alpha n+\sqrt{\alpha^2n^2+4|U||V|(1-2\alpha)}}{2}\\
&\leqslant\frac{\alpha n+\sqrt{\alpha^2n^2+4i(n-i)(1-2\alpha)}}{2}= \lambda_\alpha(K_{i,n-i}),
\end{align*}
a contradiction. Hence $G\cong K_{i,n-i}.$

Combining with Lemma~\ref{lem1.1} and \eqref{eq:3.01}, the second assertion of this theorem follows immediately. This completes the proof.
\end{proof}
Together with Theorem \ref{thm1.5} and \eqref{eq:3.01}, the next result follows immediately.
\begin{cor}
Let $G$ be a bipartite graph with order $n$ and matching number $\mu\,\left(1\leqslant \mu\leqslant \left\lfloor\frac{n}{2}\right\rfloor\right)$. If $\alpha\in[\frac{1}{2},1),$ then
$$
  \lambda_\alpha(G)\leqslant \lambda_\alpha(K_{\mu,n-\mu}).
$$
The equality holds if and only if $G\cong K_{\mu,n-\mu}.$
\end{cor}
Now, we prove the last main result of our paper.
\begin{proof}[\bf Proof of Theorem \ref{thm1.6}]
Choose an $n$-vertex graph, say $G$, with independence number $i$ such that $\lambda_\alpha(G)$ is as large as possible. Assume $I$ is a maximum independent set of $G.$ Then $|I|=i.$ By Lemma~\ref{lem3.2}, one has $G[V_G\setminus I]\cong K_{n-i}$. If we denote the graph by $G'$ which is induced by all the edges between $I$ and $V_G\setminus I$. Once again by Lemma \ref{lem3.2}, one has $G' \cong K_{i,n-i}$. Hence, $G\cong K_i^c\vee K_{n-i}$, as desired.
\end{proof}
\section{\normalsize Concluding remarks}\setcounter{equation}{0}
In this paper, we characterize the graphs with minimum $A_\alpha$-index among $\mathcal{G}_{n,i}$ for $\alpha\in[0,1)$, where $i=1,\lfloor\frac{n}{2}\rfloor,\lceil\frac{n}{2}\rceil,{\lfloor\frac{n}{2}\rfloor+1},n-3,n-2,n-1.$ We also consider the same problem among $\mathcal{G}_{n,2}$ for $\alpha\in [0,\frac{3}{4}].$  In addition, put $\alpha=\frac{1}{2}$ in Theorem \ref{thm3.2}(i) and (ii), one may obtain the corresponding results for signless Laplacian matrix, which was obtained in \cite[Theorem 3.2]{0004}. However, we cannot deduce the results for adjacency matrix (resp. signless Laplacian matrix) corresponding with  Theorem~\ref{thm3.2}(i)-(iii) (resp. Theorem \ref{thm3.2}(iii)) from our results.  Hence, it is a surprise to see that our results is non-trivial.

On the other hand, we determine the unique graph (resp. tree) on $n$ vertices with given independence number having the maximum $A_\alpha$-index with $\alpha\in[0,1)$. For the $n$-vertex bipartite graphs with given independence number, we characterize the unique graph having the maximum $A_\alpha$-index with $\alpha\in[\frac{1}{2},1).$  These results extend those of adjacency matrix \cite{Ji} and signless Laplacian matrix \cite{Sheng}.

At the end of this section, we give the following research problems:
\begin{pb}
How can we characterize the graphs with minimum $A_\alpha$-index among $\mathcal{G}_{n,i}$ for $\alpha\in[0,1)$ and $i\in [3, \lfloor\frac{n}{2}\rfloor-1]\cup [{\lfloor\frac{n}{2}\rfloor+2}, n-4]?$
\end{pb}
\begin{pb}
How can we identify the graphs with given independence number having the second maximum $A_\alpha$-index among $n$-vertex trees (resp. bipartite graphs, simple graphs) for  $\alpha\in[0,1)?$
\end{pb}
We will further the above study in the near future.
\section*{\normalsize Acknowledgement}
We take this opportunity to thank the anonymous reviewers for their critical reading of the manuscript
and suggestions which have immensely helped us in getting the article to its present form.
S.W. is financially supported by the Graduate Education Innovation Grant from Central China Normal University (Grant No. 2020CXZZ073), and the Undergraduate Innovation and Entrepreneurship Grant from Central China Normal University (Grant No. 20210409037); Y.L. is financially supported by Industry-University-Research Innovation Funding of Chinese University (Grant No. 2019ITA03033); L.S. is financially supported by the National Natural Science Foundation of China (Grant No. 11671164).

\section*{\normalsize Appendix}\setcounter{equation}{0}
\begin{proof}[\bf The proof of Claim \ref{c3}]\
By Lemma \ref{lem4.01}, it suffices to prove that the largest eigenvalue of {$A_\alpha(F(t,t+1))_{\pi_2}$} is less than $t+\frac{1}{4}.$ By a direct calculation, we have
{\begin{align*}
  \det\left(xI-A_\alpha(F(t,t+1))_{\pi_2}\right)=&x^4-(2\alpha t+2t+3\alpha-3) x^3+ ((\alpha^2+4 \alpha+1) t^2+(3 \alpha^2+ 2 \alpha- 5) t+ 3 \alpha^2\\
  &- 6 \alpha+2) x^2-(2\alpha(\alpha+1)t^3+(5 \alpha^2-3\alpha-2) t^2+(2 \alpha^3- \alpha^2-5\alpha +2)t\\
  &+\alpha^3- \alpha^2-2 \alpha+2) x+\alpha^2 t^4+ 2 \alpha(\alpha-1) t^3+ (2\alpha^3- 4 \alpha^2+ \alpha)t^2\\
  &+ (\alpha^2- 3 \alpha+2) t+2\alpha(\alpha^2- 3 \alpha+3)-2
\end{align*}
and}
\begin{align*}
  \det\left((t+\frac{1}{4})I-A_\alpha(F(t,t+1))_{\pi_2}\right)=&\frac{1}{256}\big(80(\alpha-1)^2 t^2-8\alpha(6\alpha(8\alpha-21)+103)t\\
  &+200t+ 4\alpha(4\alpha(28\alpha-89)+389)-595\big).
\end{align*}
{Let
$$
  p(x)=\det\left(xI-A_\alpha(F(t,t+1))_{\pi_2}\right)
$$
be a real function {of $x$} for $x\in [t+\frac{1}{4},+\infty),$ and let
$$
  h(t)=80(\alpha-1)^2 t^2-8\alpha(6\alpha(8\alpha-21)+103)t+200t+ 4\alpha(4\alpha(28\alpha-89)+389)-595
$$
be a real function {of $t$} for $t\in [5,+\infty).$

We firstly show that $\det\left((t+\frac{1}{4})I-A_\alpha(F(t,t+1))_{\pi_2}\right)>0.$ It is routine to check that
\begin{align*}
  &h^{(1)}(t)=160(\alpha-1)^2 t-8\alpha(6\alpha(8\alpha-21)+103)+200,\ \ \text{and}\\
  &h^{(2)}(t)=160(\alpha-1)^2>0,
\end{align*}
where $h^{(i)}({\alpha})$ refers to the $i$-th derivative function of $h(t).$ Hence $h^{(1)}(t)$ is increasing in the interval $[5,+\infty).$ Then,
$$
  h^{(1)}(t)\geqslant h^{(1)}(5)=8(1-\alpha)(2\alpha(24\alpha-89)+125)>0.
$$
So, $h(t)$ is increasing in the interval $[5,+\infty).$ Thus,
$$
  h(t)\geqslant h(5)=2405-4\alpha(4\alpha(92\alpha-351)+1641)>0.
$$
Therefore,
\[\label{eq:4.1}
\det((t+\frac{1}{4})I-A_\alpha(F(t,t+1))_{\pi_2})>0.
\]

Now, we show that $p(x)$ is a real increasing function {of $x$} for $x\in \left[t+\frac{1}{4},+\infty\right).$ It is easy to see that
\begin{align*}
p^{(1)}(x)=&4x^3-3(2\alpha t+2t+3\alpha-3) x^2+ 2((\alpha^2+4 \alpha+1) t^2+(3 \alpha^2+ 2 \alpha- 5) t+ 3 \alpha^2- 6 \alpha+2) x\\
  &-(2\alpha(\alpha+1)t^3+(5 \alpha^2-3\alpha-2) t^2+(2 \alpha^3- \alpha^2-5\alpha +2)t+\alpha^3- \alpha^2-2 \alpha+2),\\
p^{(2)}(x)=&12x^2-6(2\alpha t+2t+3\alpha-3) x+ 2((\alpha^2+4 \alpha+1) t^2+(3 \alpha^2+ 2 \alpha- 5) t+ 3 \alpha^2- 6 \alpha+2),\\
p^{(3)}(x)=&24x-6(2\alpha t+2t+3\alpha-3),\ \ \text{and}\ \ p^{(4)}(x)=24>0.
\end{align*}
Then $p^{(3)}(x)$ is an increasing function {of $x$} for $x\in[t+\frac{1}{4},+\infty).$ Thus,
$$
  p^{(3)}(x)\geqslant p^{(3)}(t+\frac{1}{4})=6(2(1-\alpha)t-3\alpha+4)>0.
$$
So, $p^{(2)}(x)$ is an increasing function {of $x$} for $x\in[t+\frac{1}{4},+\infty).$ Thus,
$$
  p^{(2)}(x)\geqslant p^{(2)}(t+\frac{1}{4})=t(2t+11)-\frac{1}{2}\alpha(2t+3)(4t+11)+2\alpha^2(t^2+3t+3)+\frac{37}{4}>0.
$$
That is, $p^{(1)}(x)$ is an increasing function {of $x$} for $x\in[t+\frac{1}{4},+\infty).$ Thus,
$$
  p^{(1)}(x)\geqslant p^{(1)}(t+\frac{1}{4})=\frac{1}{16}(24(\alpha-1)^2t^2+70t-2\alpha(4\alpha(4\alpha-17)+87)t-(5-4\alpha)^2\alpha-6)>0.
$$
That is to say, $p(x)$ is an increasing function {of $x$} for $x\in[t+\frac{1}{4},+\infty).$
Together with \eqref{eq:4.1}, we obtain that the largest eigenvalue of $A_\alpha(F(t,t+1))_{\pi_2}$ is less than $t+\frac{1}{4}.$}
This completes the proof of Claim \ref{c3}.
\end{proof}

\begin{proof}[\bf The proof of the inequality in \eqref{eq:4.4}]\
Let $h_1(\alpha)=\alpha(3+2\alpha(t-1)-2t)-3$ be a real function {of $\alpha$} for $\alpha\in[0,\frac{3}{4}{]}.$ It is routine to check that
$$
  h_1^{(1)}(\alpha)=4\alpha(t-1)-2t+3\ \text{and}\ h_1^{(2)}(\alpha)=4(t-1)>0.
$$
Hence $h_1^{(1)}(\alpha)$ is an increasing function {of $\alpha$} for $\alpha\in[0,\frac{3}{4}].$ Notice that
$$
  h_1^{(1)}(0)=-2t+3<0\ \text{and}\ h_1^{(1)}(\frac{3}{4})=t>0.
$$
So, $h_1(\alpha)$ is decreasing in the interval $\alpha\in[0,\alpha_1),$ and it is increasing in the interval $\alpha\in[\alpha_1,\frac{3}{4}],$ where $\alpha_1\in(0,\frac{3}{4}]$ with $h_1^{(1)}(\alpha_1)=0.$ Then
$$
  h_1(\alpha)\leqslant \max\{h_1(0),h_1(\frac{3}{4})\}=\max\{-3,-\frac{3}{8}(5+t)\}<0.
$$
It follows that $\det((t+1)I- A_\alpha(F(t,t+2))_{\pi})<0$ if $\alpha\in[0,\frac{3}{4}].$ This completes the proof.
\end{proof}

\begin{proof}[\bf The proof of the inequality in \eqref{eq:4.2}]\
Let $h_2(\alpha)=\det((s-1-\alpha)I- A_\alpha(F(s,t))_{\pi})$ be a real function {of $\alpha$} for $\alpha\in[0,\frac{1}{2}].$ It is routine to check that
\begin{align*}
  h_2^{(1)}(\alpha)=&\alpha\big(4\alpha^2(5+3t+s(3+t))+3\alpha(-s^2(3+t)+s(t^2-2t-11)+3(t^2+t+2))\\
    &+2(s^3+6s^2-s(t^2+4t+1)-(t-3)t-12)\big)+s(1+s)(t-s)+4s-4t+6,\\
  h_2^{(2)}(\alpha)=&\alpha\big(12\alpha(5+3t+s(3+t))+6(-s^2(3+t)+s(t^2-2t-11)+3(t^2+t+2))\big)\\
  &+2(s^3+6s^2-s(t^2+4t+1)-(t-3)t-12),\\
  h_2^{(3)}(\alpha)=&24\alpha(5+3t+s(3+t))+6(-s^2(3+t)+s(t^2-2t-11)+3(t^2+t+2)),\ \ \text{and}\\
  h_2^{(4)}(\alpha)=&24(5+3t+s(3+t))>0.
\end{align*}
Hence $h_2^{(3)}(\alpha)$ is increasing in the interval $[0,\frac{1}{2}].$ Notice that $s\geqslant t+3.$ Then
$$
  h_2^{(3)}(\alpha)\leqslant h_2^{(3)}(\frac{1}{2})=6(16-(s^2-3t)(t+3)+s(t^2-5))<0.
$$
Thus, $h_2^{(2)}(\alpha)$ is decreasing in the interval $[0,\frac{1}{2}].$ Then
$$
  h_2^{(2)}(\alpha)\geqslant h_2^{(2)}(\frac{1}{2})=2s^3-3s^2(t-1)+s(t-13)(t+2)+(t+3)(7t+3)>0.
$$
Hence $h_2^{(1)}(\alpha)$ is increasing in the interval $[0,\frac{1}{2}].$ Notice that
$$
  h_2^{(1)}(0)=s(s+1)(t-s)+4s-4t+6<0
$$
and
$$
  h_2^{(1)}(\frac{1}{2})=\frac{1}{4}(4+s^2(t+11)-s(t+1)(t+15)+t(5t+11))>0.
$$
That is to say, $h_2(\alpha)$ is decreasing in the interval $\alpha\in[0,\alpha_2),$ and it is increasing in the interval $\alpha\in[\alpha_2,\frac{1}{2}],$ where $\alpha_2\in[0,\frac{1}{2}]$ with $h_2^{(1)}(\alpha_2)=0.$ Then
$$
  h_2(\alpha)\leqslant \max\{h_2(0),h_2(\frac{1}{2})\}=\max\{t-s-1,-\frac{1}{16}(s-1)(4 s^2-6s(1+t)+t(5+2t)+1)\}<0.
$$
This completes the proof.
\end{proof}

\begin{proof}[\bf The proof of the inequality in \eqref{eq:4.3}]\
Let
$$
  h_3(\alpha)=(s-1)(1-s+t)+\alpha^2(t-1)+\alpha(s-t-2)(t-1)
$$
be a real function {of $\alpha$} for $\alpha\in[\frac{1}{2},\frac{3}{4}{]}.$ Note that $s\geqslant t+3.$ Hence
$$
  h_3^{(1)}(\alpha)=2\alpha(t-1)+(s-t-2)(t-1)>0.
$$
Hence $h_3(\alpha)$ is increasing if $\alpha\in[\frac{1}{2},\frac{3}{4}].$ Therefore,
$$
  h_3(\alpha)\leqslant h_3(\frac{3}{4})=\frac{1}{16}(-16s^2+4 s(5+7t)-t(19+12t)-1)<0.
$$
It follows that $\det((s-2+\alpha)I- A_\alpha(F(s,t))_{\pi})<0$ if $\alpha\in[\frac{1}{2},\frac{3}{4}].$ This completes the proof.
\end{proof}

\begin{thebibliography}{99}
\small \setlength{\itemsep}{-.8mm}
\bibitem{0006} M. Andeli\'{c}, C.M. Da Fonseca, S.K. Simi\'{c}, D.V. To\v{s}i\'{c}, Connected graphs of fixed order and size with maximal $Q$-index: Some spectral bounds, Discrete Appl. Math. 160 (2012) 448-459.
\bibitem{002} A. Berman, X.-D. Zhang, On the spectral radius of graphs with cut vertices, J. Combin. Theory Ser. B, 83 (2001) 233-240.
\bibitem{0009} A.E. Brouwer, W.H. Haemers, Spectra of Graphs, Springer, New York, 2012.
\bibitem{0002} R.A. Brualdi, E.S. Solheid, On the spectral radius of complementary acyclic matrices of zeros and ones, SIAM J. Algebra Discrete Method 7 (1986) 265-272.
\bibitem{001} G.X. Cai, Y. Z Fan, The signless Laplacian spectral radius of graphs with given chromatic number, Math. Appl. 22 (1) (2009) 161-167.
\bibitem{Chen} Y.Y. Chen, D. Li, J.X. Meng, On the second largest $A_{\alpha}$-eigenvalues of graphs, Linear Algebra Appl. 580 (2019) 343-358.
\bibitem{0010} E.R. van Dam, R.E. Kooij, The minimal spectral radius of graphs with a given diameter, Linear Algebra Appl. 423 (2007) 408-419.
\bibitem{RD} M. Desai, V. Rao, A characterization of the smallest eigenvalue of a graph, J. Graph Theory 18 (2) (1994) 181-194.
\bibitem{006} Y.Z. Fan, B.S. Tam, J. Zhou, Maximizing spectral radius of unoriented Laplacian matrix over bicyclic graphs of a given order, Linear and Multilinear Algebra 56 (4) (2008) 381-397.
\bibitem{007} Y.Z. Fan, D. Yang, The signless Laplacian spectral radius of graphs with given number of pendant vertices, Graphs Combin. 25 (2009) 291-298.
\bibitem{009} L.H. Feng, Q. Li, X.-D. Zhang, Spectral radii of graphs with given chromatic number, Appl. Math. Lett. 20 (2007) 158-162.
\bibitem{008} L.H. Feng, G.H. Yu, X.-D. Zhang, Spectral radius of graphs with given matching number, Linear Algebra Appl. 422 (2007) 133-138.
\bibitem{Huang} X. Huang, H.Q. Lin, J. Xue, The Nordhaus-Gaddum type inequalities of $A_{\alpha}$-matrix, Appl. Math. Comput. 365 (2020) 124716.
\bibitem{Ji} C.Y. Ji, M. Lu, On the spectral radius of trees with given independence number, Linear Algebra Appl. 488 (2016) 102-108.
\bibitem{Li} D. Li, Y.Y. Chen, J.X. Meng, The $A_{\alpha}$-spectral radius of trees and unicyclic graphs with given degree sequence, Appl. Math. Comput.  363 (2019) 124622.
\bibitem{0004} R.L. Li, L.S. Shi, The minimum signless Laplacian spectral radius of graphs with given independence number, Linear Algebra Appl. 433 (2010) 1614-1622.
\bibitem{LS3} S.C. Li, W.T. Sun, An arithmetic criterion for graphs being determined by their generalized $A_\alpha$-spectra, Discrete Math. \url{https://doi.org/10.1016/j.disc.2021.112469}.
\bibitem{LS2} S.C. Li, W.T. Sun, Some spectral inequalities for connected bipartite graphs with maximum $A_\alpha$-index, Discrete Appl. Math. 287 (2020) 97-109.
\bibitem{LS1} S.C. Li, W.T. Sun, Some bounds on the $A_\alpha$-index of connected graphs with fixed order and size,  Linear Multilinear Algebra,  \url{https://doi:10.1080/03081087.2021.1932710}.
\bibitem{LI2019}S.C. Li, S.J. Wang, The $A_{\alpha}$-spectrum of graph product, Electron. J. Linear Algebra  35  (2019) 473-481.
\bibitem{LI2020}S.C. Li, W. Wei, The multiplicity of an $A_{\alpha}$-eigenvalue: A unified approach for mixed graphs and complex unit gain graphs, Discrete Math.  343  (2020)  no. 8, 111916, 19 pp.
\bibitem{LY} S.C. Li, Y.T. Yu, The effect on $A_\alpha$-eigenvalues of mixed graphs and unit gain graphs by adding edges in clusters, Linear Multilinear Algebra,  \url{https://doi.org/10.1080/03081087.2021.1926415}.
\bibitem{0001} H.Q. Lin, X. Huang, J. Xue, A note on the $A_{\alpha}$-spectral radius of graphs, Linear Algebra Appl. 557 (2018) 430-437.
\bibitem{013} V. Nikiforov, Bounds on graph eigenvalues, II, Linear Algebra Appl. 427 (2007) 183-189.
\bibitem{0007} V. Nikiforov, Merging the $A$- and $Q$-spectral theories, Appl. Anal. Discrete Math. 11 (2017) 81-107.
\bibitem{Nik} V. Nikiforov, G. Past\'{e}n, O. Rojo, R.L. Soto, On the $A_{\alpha}$-spectra of trees, Linear Algebra Appl. 520 (2017) 286-305.
\bibitem{Ni} V. Nikiforov, O. Rojo, A note on the positive semidefiniteness of $A_{\alpha}(G)$, Linear Algebra Appl. 519 (2017) 156-163.
\bibitem{0005} V. Nikiforov, O. Rojo, On the $\alpha$-index of graphs with pendent paths, Linear Algebra Appl. 550 (2018) 87-104.
\bibitem{RJ}O. Rojo, The maximal $\alpha$-index of trees with $k$ pendent vertices and its computation,  Electron. J. Linear Algebra, 36 (2020) 38-46.
\bibitem{RJ1}Z. Lin, L.Y. Miao, S.G. Guo, Bounds on the $A_\alpha$-spread of a graph,  Electron. J. Linear Algebra, 36 (2020) 214-227.
\bibitem{Sheng} J. Sheng, M.L. Ye, The spectral radius of signless Laplacian of a connected graph with given independence number, Mathematica Applicata 23 (4) (2010) 709-712.
\bibitem{005} S.K. Simi\'{c}, F. Belardo, E.M. Li Marzi, D.V. To\v{s}i\'{c}, Connected graphs of fixed order and size with maximal index: Some spectral bounds, Linear Algebra Appl. 432 (9) (2010) 2361-2372.
\bibitem{0013} D. Stevanovi\'{c}, P. Hansen, The minimum spectral radius of graphs with a given clique number, Electron. J. Linear Algebra 17 (2008) 110-117.
\bibitem{Wang1} J.F. Wang, J. Wang, X.G. Liu, F. Belardo, Graphs whose $A_\alpha$-spectral radius does not exceed 2, Discussiones Mathematicae
Graph Theory 40 (2020) 677-690.
\bibitem{Wang} S. Wang, D. Wong, F.L. Tian, Bounds for the largest and the smallest $A_{\alpha}$ eigenvalues of a graph in terms of vertex degrees, Linear  Algebra Appl. 590 (2020) 210-223.
\bibitem{math} Wolfram Research, Inc. \textit{Mathematica}, Version 9.0 (Wolfram Research Inc., Champaign, 2012).
\bibitem{016} R.R. Wu, Y.Z. Fan, The signless Laplacian spectral radius of graphs with given number of cut edges, J. Anhui Univ. Sci. Technol. (Natural Science) 29 (2) (2009) 66-69.
\bibitem{Xu2020} F. Xu, D. Wong, F.L. Tian, On the multiplicity of $\alpha$ as an eigenvalue of the $A_{\alpha}$ matrix of a graph in terms of the number of pendant vertices, Linear Algebra Appl.  594  (2020) 193-204.
\bibitem{Xu} M.M. Xu, Y. Hong, J.L. Shu, M.Q. Zhai, The minimum spectral radius of graphs with a given independence number, Linear Algebra Appl. 431 (2009) 937-945.
\bibitem{Xue} J. Xue, H.Q. Lin, S.T. Liu, J.L. Shu, On the $A_{\alpha}$-spectral radius of a graph, Linear Algebra Appl. 550 (2018) 105-120.
\bibitem{YLH}L.H. You, M. Yang, W. So, W.G. Xi, On the spectrum of an equitable quotient matrix and its application, Linear Algebra Appl. 577 (2019) 21-40.
\end{thebibliography}
\end{document}